\pdfoutput=1
\documentclass[reqno,a4paper,11pt]{amsart}
\usepackage{amsfonts}
\usepackage{amsmath}
\usepackage{amssymb}
\usepackage[pdftitle={Distinguishing G2-manifolds}]{hyperref}
\usepackage[utf8]{inputenc}

\newcommand{\showdate}{false}

\usepackage{amssymb}
\usepackage[utf8]{inputenc}
\usepackage{a4wide}
\usepackage{enumitem}
\usepackage{tikz}
\usepackage{multido}
\usepackage{color}
\usepackage{microtype}
\usepackage{ifthen}
\usepackage{comment}
\usepackage[all]{xy}
\RequirePackage{xspace}

\ifthenelse{\boolean{\showdate}}{\date{\today}}{}

\newcommand{\ignore}[1]{}

\DeclareMathOperator{\ind}{ind}
\DeclareMathOperator{\im}{Im}
\DeclareMathOperator{\rk}{rk}
\DeclareMathOperator{\lcm}{lcm}
\DeclareMathOperator\tr{tr}
\DeclareMathOperator\arc{arc}
\renewcommand\Re{\operatorname{Re}}

\newcommand\Spin{\mathrm{Spin}}

\newcommand{\ie}{\emph{i.e.} }
\newcommand{\eg}{\emph{e.g.} }
\newcommand\cf{see }
\newcommand{\signat}{\sigma}
\newcommand{\kd}{\Sigma}

\newcommand{\hk}{hyper-Kähler\xspace}
\DeclareMathAlphabet{\matheur}{U}{eur}{m}{n}
\newcommand{\tormat}{\matheur{t}}
\newcommand{\hkr}{\matheur{r}}

\newcommand{\PP}{\mathbb{P}}

\newcommand{\ahat}{\widehat{A}}
\newcommand{\mmod}{\!\!\mod}

\newcommand{\C}{\mathbb{C}}
\newcommand{\Z}{\mathbb{Z}}
\newcommand{\Q}{\mathbb{Q}}

\newcommand{\R}{\mathbb{R}}

\newcommand{\bbrp}{\mathbb{R}^{+}}

\newcommand{\into}{\hookrightarrow}

\newcommand{\Num}[1]{\textup{Num}\left(#1\right)}

\newcommand{\del}{\partial}

\newcommand{\gstr}{$G$\nobreakdash-\hspace{0pt}structure}
\newcommand{\gtstr}{$G_{2}$\nobreakdash-\hspace{0pt}structure}

\newcommand{\gtmfd}{$G_{2}$\nobreakdash-\hspace{0pt}manifold}

\newcommand{\gtmetric}{$G_2$-metric}

\newcommand{\calz}{\mathcal{Z}}
\newcommand{\co}{\mathcal{O}}

\newcommand{\anglen}{u}
\newcommand{\anglex}{v}
\newcommand{\lnx}{\xi}
\newcommand{\lnn}{\zeta}

\newcommand\sign{\operatorname{sign}}
\newcommand\abs[1]{\left|#1\right|}

\newcommand\Sign{\operatorname{sign}}

\newcommand\gmatrix[4]{\bigl(\begin{smallmatrix}#1&#2\\#3&#4\end{smallmatrix}\bigr)}

\let\eps\varepsilon
\let\phy\varphi
\let\thet\vartheta
\let\<\langle
\let\>\rangle

\newcommand\K{\Sigma} %

\newtheorem{thm}{Theorem}[section]

\newtheorem{cor}[thm]{Corollary}

\theoremstyle{definition}
\newtheorem{dfn}[thm]{Definition}

\theoremstyle{remark}
\newtheorem{rmk}[thm]{Remark}
\newtheorem{ex}[thm]{Example}

\setlist{leftmargin=*}
\newcommand{\td}{\mkern3mu{\tilde {\mkern-3mu d \mkern-1mu}}\mkern1mu}

\definecolor{darkgreen}{rgb}{0,0.5,0}
\definecolor{darkred}{rgb}{0.7,0,0}
\definecolor{darkblue}{rgb}{0,.2,.7}

\newcommand{\thefigures}{
\protect \begin{figure}
\protect \begin{minipage}{0.48\textwidth}
\centering
\begin{tikzpicture}[x=0.8cm,y=0.8cm]
  \multido{\ix=-2+2}{3}{
    \multido{\iy=-2+2}{3}{
      \fill (\ix,\iy) circle (2pt) ;
    }
  }
  \multido{\ix=-1+2}{2}{
    \multido{\iy=-1+2}{2}{
      \fill (\ix,\iy) circle (2pt) ;
    }
  }
  \draw (0,0) -- (0,2) ;
  \draw (0,0) -- (1,-1) ;
  \draw (0,0) -- (1,1) ;
  \draw (0,0) -- (2,0) ;
  \begin{scope}[->, line width=1pt]
    \draw[color=blue] (0,0) -- (90:1.131cm)
	node[above left, color=black] {$\partial_{u_+}$} ;
    \draw[color=blue] (0,0) -- (-45:1.131cm)
	node[below right, color=black] {$\partial_{u_-}$} ;
    \draw[color=red] (0,0) -- (45:1.131cm)
	node[above right, color=black] {$\partial_{v_-}$} ;
    \draw[color=red] (0,0) -- (0:1.131cm)
	node[below right, color=black] {$\partial_{v_+}$} ;
  \end{scope}
  \fill (0,0) circle (3pt) ;
  \draw[->] (0,0) ++(0:0.6cm) arc (0:45:0.6cm) ;
  \node at (22.5:0.85cm) {$\vartheta$} ;
\end{tikzpicture}
\caption{\texorpdfstring{$G=\protect\gmatrix111{-1},\vartheta = \tfrac{\pi}{4}$}{angle pi/4}}
\label{fig:1/4}
\protect \end{minipage}\hfill
\protect \begin{minipage}{0.48\textwidth}
\centering
\begin{tikzpicture}[x=1cm,y=0.577cm]
  \multido{\ix=-2+2}{3}{
    \multido{\iy=-2+2}{3}{
      \fill (\ix,\iy) circle (2pt) ;
    }
  }
  \multido{\ix=-1+2}{2}{
    \multido{\iy=-3+2}{4}{
      \fill (\ix,\iy) circle (2pt) ;
    }
  }
  \draw (0,0) -- (0,2) ;
  \draw (0,0) -- (1,-3) ;
  \draw (0,0) -- (1,1) ;
  \draw (0,0) -- (2,0) ;
  \begin{scope}[->, line width=1pt]
    \draw[color=blue] (0,0) -- (90:1.154cm)
	node[above, color=black] {$\partial_{u_+}$} ;
    \draw[color=blue] (0,0) -- (-60:1.154cm)
	node[right, color=black] {$\partial_{u_-}$} ;
    \draw[color=red] (0,0) -- (30:1.154cm)
	node[above right, color=black] {$\partial_{v_-}$} ;
    \draw[color=red] (0,0) -- (0:1.154cm)
	node[below right, color=black] {$\partial_{v_+}$} ;
  \end{scope}
  \fill (0,0) circle (3pt) ;
  \draw[->] (0,0) ++(0:0.75cm) arc (0:30:0.75cm) ;
  \node at (15:1cm) {$\vartheta$} ;
\end{tikzpicture}
\caption{\texorpdfstring{$G=\protect\gmatrix111{-3},\vartheta = \tfrac{\pi}{6}$}{angle pi/6}}
\label{fig:1/6}
\protect \end{minipage}
\protect \end{figure}
}

\newcommand{\morefigures}{
\protect \begin{figure}
\protect \begin{minipage}{0.48\textwidth}
\centering
\begin{tikzpicture}[x=1cm,y=0.577cm]
  \multido{\ix=-2+2}{3}{
    \multido{\iy=-2+2}{3}{
      \fill (\ix,\iy) circle (2pt) ;
    }
  }
  \multido{\ix=-1+2}{2}{
    \multido{\iy=-3+2}{4}{
      \fill (\ix,\iy) circle (2pt) ;
    }
  }
  \draw (0,0) -- (0,2) ;
  \draw (0,0) -- (1,-1) ;
  \draw (0,0) -- (1,3) ;
  \draw (0,0) -- (2,0) ;
  \begin{scope}[->, line width=1pt]
    \draw[color=blue] (0,0) -- (90:1.154cm)
	node[above, color=black] {$\partial_{u_+}$} ;
    \draw[color=blue] (0,0) -- (-30:1.154cm)
	node[below right, color=black] {$\partial_{u_-}$} ;
    \draw[color=red] (0,0) -- (60:1.154cm)
	node[right, color=black] {$\partial_{v_-}$} ;
    \draw[color=red] (0,0) -- (0:1.154cm)
	node[below right, color=black] {$\partial_{v_+}$} ;
  \end{scope}
  \fill (0,0) circle (3pt) ;
  \draw[->] (0,0) ++(0:0.6cm) arc (0:60:0.6cm) ;
  \node at (30:0.85cm) {$\vartheta$} ;
\end{tikzpicture}
\caption{\texorpdfstring{$G=\protect\gmatrix113{-1},\vartheta = \tfrac{\pi}{3}$}{angle pi/3}}
\label{fig:1/3}
\protect \end{minipage}\hfill
\protect \begin{minipage}{0.48\textwidth}
\centering
\begin{tikzpicture}[x=0.816cm,y=0.577cm]
  \multido{\ix=-3+3}{3}{
    \multido{\iy=-3+3}{3}{
      \fill (\ix,\iy) circle (2pt) ;
    }
  }
  \multido{\ix=-2+3}{2}{
    \multido{\iy=-1+3}{2}{
      \fill (\ix,\iy) circle (2pt) ;
    }
  }
  \multido{\ix=-1+3}{2}{
    \multido{\iy=-2+3}{2}{
      \fill (\ix,\iy) circle (2pt) ;
    }
  }
  \draw (0,0) -- (0,3) ;
  \draw (0,0) -- (1,-1) ;
  \draw (0,0) -- (1,2) ;
  \draw (0,0) -- (3,0) ;
  \begin{scope}[->, line width=1pt]
    \draw[color=blue] (0,0) -- (90:1cm)
	node[above left, color=black] {$\partial_{u_+}$} ;
    \draw[color=blue] (0,0) -- (-35.3:1cm)
	node[below right, color=black] {$\partial_{u_-}$} ;
    \draw[color=red] (0,0) -- (54.7:1cm)
	node[right, color=black] {$\partial_{v_-}$} ;
    \draw[color=red] (0,0) -- (0:1cm)
	node[below right, color=black] {$\partial_{v_+}$} ;
  \end{scope}
  \fill (0,0) circle (3pt) ;
  \draw[->] (0,0) ++(0:0.6cm) arc (0:54.7:0.6cm) ;
  \node at (27.3:0.85cm) {$\vartheta$} ;
\end{tikzpicture}
\caption{\texorpdfstring{$G=\protect\gmatrix112{-1},\vartheta = \arc\cos \tfrac1{\sqrt 3}$}{angle pi/6}}
\label{fig:A4}
\protect \end{minipage}
\protect \end{figure}
}

\newcommand\hypfigure{
  \protect\begin{figure}
      \begin{tikzpicture}
    \begin{scope}
      \clip  (-4,0) arc (180:70.5:3) -- (0,0)
      arc (0:180:0.5)
      arc (0:180:1.5) -- cycle ;
      \draw (-1,0) ++(70.5:3) circle (0.7) ;
    \end{scope}
    \draw (-1,0) ++(70.5:3) ++(215:0.35) node {$2\thet$} ;
    \draw (2,0) node[below] {$0$} arc (0:70.5:3) ;
    \draw[->] (0,0) -- (0,3.5) node[above] {$\infty$} ;
    \draw[color=green, line width=1.5pt] (-4,0)
      arc (180:70.5:3) node[pos=0.5,above,color=black] {$\gamma_-$} --
      node[pos=0.5,right,color=black] {$\gamma_+$} (0,0) ;
    \draw[color=green] (0,0) node[below,color=black] {$-\frac13$}
      arc (0:180:0.5) node[below,color=black] {$-\frac12$}
      arc (0:180:1.5) node[below,color=black] {$-1$} ;
    \draw (-5,0) -- (3,0) node[below right] {$\del_\infty\mathcal H$} ;
    \node at (-1,1.65) {$P$} ;
      \end{tikzpicture}
    \caption{The hyperbolic polygon for Example~\protect\ref{ex:hyp}}\label{fig:hyp}
  \end{figure}
}

\begin{document}

\title{Distinguishing $G_2$-manifolds}

\author{Diarmuid Crowley, Sebastian Goette and Johannes Nordstr\"om}


\address{School of Mathematics and Statistics,
University of Melbourne,
Parkville, VIC, 3010, \mbox{Australia}}
\email{dcrowley@unimelb.edu.au}

\address{Mathematisches Institut\\ Universit\"at Freiburg\\ Eckerstr.~1, 79104 Freiburg, Germany} \email{sebastian.goette@math.uni-freiburg.de}

\address{Department of Mathematical Sciences\\
University of Bath\\
Bath BA2 7AY\\
UK}\email{j.nordstrom@bath.ac.uk}

\begin{abstract}
In this survey, we describe invariants that can be used to distinguish
connected components of the moduli space of holonomy $G_2$ metrics on a closed
7-manifold, or to distinguish \gtmfd s that are homeomorphic but not
diffeomorphic.  We also describe the twisted connected sum and extra-twisted
connected sum constructions used to realise $G_2$-manifolds for which the above
invariants differ.
\end{abstract}

\maketitle

\vspace{-0.75\baselineskip}

\section{Introduction}

This is a survey of recent results on the topology of closed Riemannian
7-manifolds with holonomy $G_2$ and their \gtstr s.
Among the highlights are examples of
\begin{itemize}
\item closed 7-manifolds whose moduli space of holonomy $G_2$ metrics is
disconnected, \ie the manifold admits a pair of \gtmetric s that cannot
be connected by a path of \gtmetric s (even after applying a diffeomorphism
to one of them)
\item pairs of closed 7-manifolds that both admit holonomy $G_2$ metrics,
which are homeomorphic but not diffeomorphic
\end{itemize}
The key ingredients are
\begin{itemize}
\item invariants that can distinguish homeomorphic closed 7-manifolds up to
diffeomorphism, or \gtstr s on 7-manifolds up to homotopy and diffeomorphism
\item classification theorems for smooth 7-manifolds or \gtstr s on smooth 7-manifolds, at least for 7-manifolds that are 2-connected
\item a method for producing many examples of closed \gtmfd s, many of which
are 2-connected, and for which the above invariants can be computed
\end{itemize}
By a homotopy of \gtstr s we simply mean a continuous path of $G_2$-structures
on a fixed manifold. The relevance is that metrics with holonomy $G_2$ are
essentially equivalent to \emph{torsion-free} \gtstr s.
If two \gtmetric s on $M$
are in the same component of the moduli space, then their associated \gtstr s
are certainly related by homotopy and diffeomorphism. However, studying
homotopy classes of \gtstr s is essentially a topological problem, and avoids
considering the complicated partial differential equation of torsion-freeness.

This survey will concentrate on describing the invariants and the
constructions, while only stating the classification results.

\enlargethispage{\baselineskip}

\subsection{Invariants and classification results for 2-connected 7-manifolds}

By Poincar\'e duality, all information about the cohomology of a closed
2-connected 7-manifold $M$ is captured by $H^4(M)$. For simplicity, let us from
now on assume that $H^4(M)$ is torsion-free. Then in particular the data about
the cohomology reduces to the integer $b_3(M)$.

A 2-connected manifold has a unique spin structure. The only interesting
relevant characteristic class of $M$ is the spin characteristic class
$p_M \in H^4(M)$
(which determines the first Pontrjagin class by $p_1(M) = 2p_M$).
Since we assume $H^4(M)$ to be torsion-free, the data of $p_M$ amounts
to specifying the greatest integer $d$ that divides $p_M$ in $H^4(M)$
(where we set $d := 0$ if $p_M = 0$).
This $d$ is in fact even, see \S\ref{subsec:p}.

\begin{thm}[{\cite[Theorem 3]{wilkens71}}]
Closed 2-connected 7-manifolds with torsion-free cohomology are classified up
to homeomorphism by $(b_3, d)$.
\end{thm}

If $p_M = 0$, %
or more generally, if~$p_M$ is a torsion class, %
then $M$ admits 28 diffeomorphisms classes of smooth structures,
distinguished by the diffeomorphism invariant of Eells and
Kuiper \cite{eells62}.
A generalisation of this invariant to the case when $p_M$ is non-torsion was
introduced in \cite{7class}. Under the simplifying assumption that $H^4(M)$ is
torsion-free, this generalised Eells-Kuiper invariant is a constant
\[ \mu(M) \in \Z/\gcd(28, \tfrac{\td}{4})\;, \]
where
\[ \td := \lcm(4, d)\; . \]

\begin{thm}[{\cite[Theorem 1.3]{7class}}]
Closed 2-connected 7-manifolds with torsion-free cohomology are classified up
to diffeomorphism by $(b_3, d, \mu)$.
\end{thm}

In particular, the number of diffeomorphism classes of smooth structures on
a 2-connected $M$ is exactly $\gcd(28, \tfrac{\td}{4})$.

Given a \gtstr{} on $M$, \cite{nu} defines %
two further %
invariants $\nu$
and $\xi$. The first is simply a constant $\nu \in \Z/48$. Under the
assumption that $H^4(M)$ is torsion-free, $\xi$ is also a constant
$\xi(M) \in \Z/3\td$.
Both are invariant not only under %
diffeomorphisms %
but also under homotopies of
\gtstr s. %
They satisfy the relations
\begin{subequations}
\label{eq:constraints}
\begin{align}
\label{eq:parity}
\nu &= \sum_{i = 0}^3 b_i(M) \mod 2 \\
\label{eq:mu_constraint}
12\mu &= \xi-7\nu \mod \gcd(12\cdot28, 3\td)\; . 
\end{align}
\end{subequations}

\begin{thm}[{\cite[Theorem 6.9]{nu}}]
\label{thm:g2class}
Closed 2-connected 7-manifolds with torsion-free cohomology equipped with a
\gtstr{} are classified up to diffeomorphism and homotopy
by $(b_3, d, \nu, \xi)$.
\end{thm}

In particular, the number of classes of \gtstr s modulo homotopy and
diffeomorphism on a fixed smooth 2-connected $M$ is determined by computing the
number of pairs $(\nu, \xi)$ that satisfy \eqref{eq:constraints}; for each of
the 24 values of $\nu$ allowed by the parity constraint \eqref{eq:parity},
there are $\Num{\frac{d}{112}}$ values for $\xi$ that satisfy
\eqref{eq:mu_constraint}, so the number of classes
is $24 \, \Num{\frac{d}{112}}$. %
Here, ``Num'' denotes the numerator of a fraction written in lowest terms. %
We can say that $\nu$ on its own always
distinguishes at least 24 classes of \gtstr s on any fixed $M$, and if $d$
divides 112 then it determines the classes completely.

The invariants $\mu$, $\nu$ and $\xi$ are all defined as ``coboundary defects''
of characteristic class formulas valid for closed 8-manifolds. The definitions
of $\nu$ and $\xi$ rely on interpreting \gtstr s in terms of non-vanishing
spinor fields.

The $\nu$-invariant is a bit more robust than the other two, in that its range
does not depend on~$d$.
If the \gtstr\ is torsion-free,
it is also possible to define a closely related
invariant $\bar \nu \in \Z$ in terms of spectral invariants of the metric
induced by the \gtstr{}, %
which satisfies
\[  \nu = \bar\nu + 24\,(1+b_1(M)) \mod 48\;, \]
see Corollary~\ref{cor:nubar}.
The analytic refinement is invariant under diffeomorphisms, but \emph{not}
under arbitrary homotopies of \gtstr s. However, $\bar \nu$ \emph{is} invariant
under homotopies through torsion-free \gtstr s. Therefore $\bar \nu$ is
capable of distinguishing components of the moduli space of $G_2$-metrics
on a manifold $M$, even when the associated \gtstr s are homotopic.

\subsection{Twisted connected sums}\label{Sect:tcs}

The source of examples that we use is the twisted connected sum construction
pioneered by Kovalev \cite{kovalev03} and studied further in \cite{g2m},
and the ``extra-twisted'' generalisation from \cite{xtcs,eta}.
Let us first outline the original version of the construction.

Suppose that $V_+$ and $V_-$ is a pair of asymptotically cylindrical Calabi-Yau
3-folds: Ricci-flat Kähler 3-folds with an asymptotic end exponentially close to
a product cylinder $\bbrp \times U$.
We require the asymptotic cross-section $U_\pm$
of $V_\pm$ to be of the form $S^1 \times \kd_\pm$ where $\kd_\pm$ is a
K3 surface. Then $M_\pm := S^1 \times V_\pm$ is an %
asymptotically cylindrical (ACyl) %
\gtmfd{} with asymptotic cross-section $Y_\pm=T_\pm \times \kd_\pm$,
where the 2-torus $T_\pm$ is a product
of an `internal' circle factor from the asymptotic cross-section of $V_\pm$
and the `external' circle factor in the definition of $M_\pm$.
Let $\tormat : T_+ \to T_-$ be an orientation-reversing isometry that
swaps the internal and external circle directions.
We call $\hkr : \kd_+ \to \kd_-$ a \emph{\hk rotation} if the product map
\begin{equation}
\label{eq:cylmap}
(-1) \times \tormat \times \hkr :
\R \times T_+ \times \kd_+ \to \R \times T_- \times \kd_-
\end{equation}
is an isomorphism of the asymptotic limits of the torsion-free \gtstr s
of $M_+$ and $M_-$ (see Definition \ref{def:hkr}).
Given a \hk rotation $\hkr$ and a sufficiently
 large `neck length' parameter $\ell$, we can truncate the cylinders of $M_\pm$
at distance $\ell$, form a closed 7-manifold $M_\ell$ by gluing the
boundaries using $\tormat \times \hkr$ and patch the torsion-free \gtstr s %
from the halves to a %
closed %
\gtstr~$\phy$ with small torsion on $M_\ell$.
By Kovalev \cite[Theorem 5.34]{kovalev03} or the more general results of %
Joyce~\cite[Theorem 11.6.1]{joyce00},
$\phy_\ell$ can be perturbed to a torsion-free \gtstr~$\bar\phy$.

The cohomology of the twisted connected sum $M_\ell$ can be computed from
the cohomology of $V_+$ and $V_-$ using Mayer-Vietoris, given some data about
the action of $\hkr$ on cohomology. It is convenient to describe the latter
piece of data in terms of what we call the \emph{configuration} of $\hkr$.
Call the image $N_\pm \subset H^2(\kd_\pm;\Z)$ of the restriction map
$H^2(V_\pm;\Z) \to H^2(\kd_\pm;\Z)$ the \emph{polarising lattice} of $V_\pm$.
If $L$ is an even unimodular lattice of signature $(3,19)$, then
$H^2(\kd_\pm;\Z) \cong L$ by the classification of lattices, so we can
identify $N_+$ and $N_-$ with sublattices of $L$, each well-defined up
to the action of the isometry group $O(L)$. Given $\hkr$ we can instead
consider the pair of embeddings $N_+, N_- \into L$ as well-defined
up to the action of $O(L)$, and we call such a pair a configuration of the
polarising lattices.

In a similar way, $p_M$ can be computed from data about $V_+$ and $V_-$
together with the configuration. %
But even without considering that data, there
are some strong general restrictions on the possible values of the greatest
divisor $d$ of $p_M$. It follows from \cite[Proposition 10.2.7]{joyce00}
that $p_M$ is rationally non-trivial, so $d > 0$. %
Note also that $M$ always
contains a K3 surface with trivial normal bundle.
Since $p_{K3} \in H^4(K3;\Z) \cong \Z$ corresponds
to 24, $d$ must always divide 24. %
As explained in section~\ref{subsec:p}, $d$ is always even, %
so \emph{a priori} the only possible values for $d$ are $2,4,6, 8, 12$ and $24$.
 
As explained in \S\ref{subsec:match}, it is easier to find
examples of pairs of ACyl Calabi-Yau 3-folds $V_+, V_-$ with a \hk rotation
of their asymptotic K3s where the configuration is `perpendicular' (in the sense
that every element of $N_+$ is perpendicular to every element of $N_-$ in $L$)
than where it is not. In~\cite{g2m} it is shown that there are at least $10^8$ pairs
$V_+, V_-$ with a perpendicular matching, which are 2-connected with $H^4(M)$
torsion-free. Computing $b_3$ and $d$ shows that many of the resulting twisted
connected sums are homeomorphic.
However, these examples on their own turn out to be insufficient for 
addressing the questions above.

If $M$ is 2-connected with torsion-free $H^4(M)$ and
$d = 2, 4, 6$ or $12$ then $M$ has a unique smooth structure (up to
diffeomorphism), while if $d = 8$ or $24$ then $M$ admits exactly two classes
of smooth structure distinguished by the generalised Eells-Kuiper invariant
$\mu \in \Z/2$. For twisted connected sums, $\mu$ is
computed in \cite{exotic} in terms of the same data used to determine~$p_M$.
It turns out that $\mu = 0$ for any twisted connected sum with
perpendicular configuration.
However, \cite{exotic} studies the problem of finding twisted connected
sums with non-perpendicular configuration, and thereby also produces some
examples with $\mu = 1$.

\begin{ex}
\label{ex:exotic}
The smooth 2-connected 7-manifolds with torsion-free $H^4(M)$ and
$(b_3, d, \mu) = (101, 8, 0)$ and $(101, 8, 1)$ both admit metrics with
holonomy $G_2$ (see Example \ref{ex:twistedcubic}); they form a pair of
\gtmfd s that are homeomorphic but not diffeomorphic.
\end{ex}

Turning to the $G_2$-moduli space of twisted connected sums,
we find that if one attempts to distinguish components
of the moduli space using the $\nu$-invariant,
the $\nu$-invariant of a twisted connected sum turns out to always take the
same value.
Indeed, it was computed in \cite{nu} in terms of a spin cobordism that all
twisted connected sums have $\nu = 24$, regardless of the ACyl Calabi-Yaus used
or the configuration.
Considering the analytic refinement $\bar\nu$ does not help either.

\begin{thm}[{\cite[Corollary~3]{eta}}] \label{thm:nu=0}
$\bar\nu = 0$ for any twisted connected sum.
\end{thm}

If we want to use this circle of ideas to exhibit examples of closed
7-manifolds with disconnected moduli space of $G_2$ metrics, we are left with
two possible approaches. One is to make use of the $\xi$-invariant, and
this approach has 
very recently been successfully followed by Wallis
\cite{wallis18}. His computation of the $\xi$-invariant for twisted connected
sums shows that, %
like~$\mu$, %
it is uninteresting whenever the configuration is
perpendicular. In particular, none of the examples found in \cite{g2m} can be
distinguished using $\xi$. However, \cite[Examples 1.6~\&~1.7]{wallis18}
provide twisted connected sums with non-perpendicular configuration 
and $d = 6$ or~$24$, where $\xi$ does distinguish components of the moduli space. %

The other approach available for disconnecting the $G_2$-moduli space
is to consider a more general class of examples and we review recent
work along this direction in the following subsection.
Among these new examples we will also find $G_2$-manifolds
that are not $G_2$-nullbordant,
see Remark~\ref{rem:G2bordism} below,
so $G_2$-bordism presents no obstruction against holonomy~$G_2$.

\subsection{Extra-twisted connected sums}\label{xtcsSect}

Our generalisation of the twisted connected sum construction relies on
using ACyl Calabi-Yau manifolds $V_\pm$ with automorphism groups
$\Gamma_\pm \cong \Z/k_\pm$. The action of $\Gamma_\pm$ on the asymptotic
cross-section $S^1 \times \kd_\pm$ of $V_\pm$ is required to be trivial
on $\kd_\pm$ and free on the `internal' $S^1$ factor.
If we also let $\Gamma_\pm$ act freely on the `external' $S^1$ factor of
$S^1 \times V_\pm$, then the quotient $M_\pm := (S^1 \times V_\pm)/\Gamma_\pm$
is a smooth ACyl \gtmfd. The asymptotic cross-section is of the form
$T_\pm \times \kd_\pm$, for $T_\pm := (S^1 \times S^1)/\Gamma_\pm$.
Note that the 2-torus $T_\pm$ need not be a metric product of two circles;
on the other hand it could be even if $k_\pm > 1$, depending on the choice
of circumferences of the internal and external circles.

Suppose we have arranged the circumferences of the circles in such a
way that there exists an orientation-reversing isometry
$\tormat : T_+ \to T_-$. A key parameter of $\tormat$ is the angle $\thet$
between the external circle directions. For a diffeomorphism
$\hkr : \kd_+ \to \kd_-$, the condition that \eqref{eq:cylmap} be an
isomorphism of the asymptotic limits of $M_+$ and $M_-$ depends on $\thet$
(see Definition \ref{def:hkr}). Given such a $\thet$-\hk rotation, we can
proceed to glue $M_+$ and $M_-$ similarly %
as before to form a closed manifold
$M$ with a torsion-free \gtstr.
We assume that $\thet$ is not a multiple of $\pi$, so that $M$ has finite
fundamental group, and thus holonomy exactly $G_2$
(otherwise the external circles are aligned, so that $M$ has an $S^1$ factor).
Note that if $k_+ = k_- = 1$ then~$\thet$ is forced to be a right angle
and we recover the ordinary twisted connected sum construction from the
previous subsection.

\enlargethispage{2\baselineskip}

Unlike for ordinary twisted connected sums, the analytic invariant of
an extra-twisted connected sum \emph{is} affected by the configuration.
In fact, when both $k_\pm \leq 2$, the only contributions to $\bar \nu$
come from~$\rho=\pi-2\thet$, and
the invariant $m_\rho(L;N_+,N_-)\in\Z$ of the configuration defined in
\eqref{eq:mrho} (see \cite[Def~2.5]{eta}).

If~$k_\pm\ge 3$, there are two more contributions.
The generalised Dedekind sum~$D_{\gamma_\pm}(V_\pm)\in\Q$
defined in \eqref{eq:dedekind} (see~\cite%
{nuxx})
depends on the action of~$\Gamma_\pm$ on~$V_\pm$.
It vanishes if no element~$\gamma\in\Gamma_\pm$ has isolated %
fixed points. %
On an odd-dimensional Calabi-Yau manifold,
no structure preserving involution can have isolated %
fixed points, %
so this contribution vanishes if~$k_\pm\le 2$.

Finally, \eqref{eq:F} defines a number~$F_\pm\in\R$
that depends on the circumferences of the internal and external~$S^1$
and the $\Gamma_\pm$-action on their product
(see~\cite%
{nuxx}).
It vanishes in the rectangular case ($k_\pm=1$)
and in the rhombic case ($k_\pm=2$).
While it is hard to compute $F_\pm$ individually in general,
we sketch ways to determine~$F_+ + F_-$ in Sections~\ref{sect:nuxxtcs}
and~\ref{sec:hyp} below.

\begin{thm}[{\cite[Thm~1]{eta}, see also
      \cite%
          {nuxx}}]\label{thm:involutions}
  Let~$(M,g)$ be an extra-twisted connected sum. %
  Let $\thet$ be the gluing angle, and let~$m_\rho(L;N_+,N_-)$,
  $D_{\gamma_\pm}(V_\pm)$ and~$F_\pm$ be as above.
  Then
  \begin{equation*}
      \bar\nu(M)
      =D_{\gamma_+}(V_+)+D_{\gamma_-}(V_-)
      +F_+ + F_-
      -72\frac\rho\pi+3m_\rho(L;N_+,N_-)\;.
  \end{equation*}
\end{thm}

\begin{ex}\label{ex36}
  The smooth 2-connected 7-manifold with torsion-free~$H^4(M)$
  and~$(b_3,d,\mu)=(97,2,0)$ admits two torsion-free $G_2$-structures
  with~$\bar\nu=0$ and~$\bar\nu=-36$, see~\cite[Ex~3.7]{eta}.
  Hence these torsion-free \gtstr s are not homotopic, so the corresponding
  holonomy $G_2$ metrics must lie in different components of the $G_2$ moduli
  space.
  One of the two $G_2$-structures comes from an extra-twisted connected %
  sum with gluing angle~$\thet = \frac\pi 4$, while the other is a rectangular
  twisted connected sum.
\end{ex}

\begin{ex}\label{ex48}
  The smooth 2-connected 7-manifold with torsion-free~$H^4(M)$
  and~$(b_3,d,\mu)=(109,2,0)$ admits two torsion-free $G_2$-structures
  with~$\bar\nu=0$ and~$\bar\nu=-48$, see~\cite[Ex~3.11]{eta}.
  Both have~$\nu=24$, and because~$d$ divides~$112$, the underlying
  $G_2$-structures are homotopic (after choosing the diffeomorphism
  appropriately) by Theorem \ref{thm:g2class}.
  Nevertheless, the analytic invariant~$\bar\nu$ %
  shows that the corresponding holonomy
  $G_2$ metrics are in different components of the $G_2$ moduli space.
  One of the two $G_2$-structures comes from an extra-twisted connected %
  sum with gluing angle~$\frac\pi 6$, while the other is a rectangular twisted
  connected sum.
\end{ex}

\begin{rmk}
  The examples above all have~$d=\td=4$ and~$12|\nu$,
  so they have~$\xi=0$ by~\eqref{eq:mu_constraint}.
  On the other hand, the examples found by Wallis all have $\bar \nu = 0$
  by Theorem \ref{thm:nu=0}, demonstrating that neither~$\xi$ nor~$\bar \nu$
  is a complete invariant of the connected components of the 
  $G_2$-moduli space. (Nor do we have any reason to believe that these methods
  can give a complete set of invariants %
  of the connected components.) 
\end{rmk}

\begin{rmk}\label{rem:G2bordism}
  The $\nu$-invariants in Theorem~\ref{thm:involutions} are always
  divisible by~$3$ if~$k_\pm\le 2$.
  If~$k_\pm>2$, it is possible to construct examples where this is no longer
  the case, see Example~\ref{ex:xbord}.
  Indeed, we %
  expect %
  that~$\nu$ can attain all values in~$\Z/48$
  satisfying the parity constraint~\eqref{eq:parity}.

  This is significant because a topological $G_2$-structure is trivial
  in $G_2$-bordism if and only if~$3|\nu$ (or equivalently, $3|\xi$),
  see Remark \ref{rmk:g2cob}. 
  Hence we see that $G_2$-bordism does not give an obstruction against
  the existence of torsion-free $G_2$-structures.
\end{rmk}

\enlargethispage{\baselineskip}

\begin{rmk}
  The fact that~$\bar\nu(M)$ is always an integer poses interesting
  restrictions on the possible asymptotically cylindrical pieces~$V_\pm$,
  the groups~$\Gamma_\pm$,
  the torus matchings, and the matchings of K3 surfaces.
  For~$k_\pm\ge 3$, the values of~$F_\pm$ and~$\frac\rho\pi$
  can be irrational (see Figure~\ref{fig:A4}).
  Using elementary hyperbolic geometry, one can prove that the linear
  combination of these terms that occurs in Theorem~\ref{thm:involutions}
  is always rational.

  The generalised Dedekind sums~$D_{\gamma_\pm}(V_\pm)$ are always rational, too.
  The fractional part of their sum is determined by the remaining
  terms in Theorem~\ref{thm:involutions}.
  As an example, if an asymptotically cylindrical Calabi-Yau manifold~$V_\pm$
  with an action of~$\Gamma_\pm\cong\Z/5\Z$ occurs in a matching,
  then the action of~$\Gamma_\pm$ on that space must have isolated %
  fixed points. %
\end{rmk}

\begin{rmk}
  Finally, one may wonder if it is worthwhile to extend the construction above
  by allowing group actions of~$\Gamma_\pm$ on~$V_\pm$ that do not necessarily
  act trivially on the K3-factor~$\Sigma_\pm$.
  While on one hand it may be difficult to provide such examples,
  it turns out that on the other hand this will only give quotients
  of the examples we can produce by our methods above.

  To understand this, let~$\Gamma_{0,\pm}\subset\Gamma_{\pm}$
  be the normal subgroup of~$\Gamma_\pm$ that fixes~$\Sigma_\pm$ pointwise.
  The cross-section of~$(V_\pm\times S^1)/\Gamma_\pm$
  at infinity can be regarded as %
  the %
  total space of a singular fibration
  \begin{equation*}
    (S^1\times S^1\times\K_\pm)/\Gamma_\pm\longrightarrow
    \K_\pm/(\Gamma_\pm/\Gamma_{0,\pm})
  \end{equation*}
  that is locally of product geometry, and whose regular fibres are
  all isometric to~$(S^1\times S^1)/\Gamma_{0,\pm}$.
  
  Now assume that we can glue~$(S^1\times V_-)/\Gamma_-$
  to~$(S^1\times V_+)/\Gamma_+$, obtaining a $G_2$-manifold~$M$.
  Then the isometry of the cross-sections at infinity on both sides lifts
  to an isometry
  \begin{equation*}
    (S^1\times S^1)/\Gamma_{0,-}\times K_-\stackrel\cong\longrightarrow
    (S^1\times S^1)/\Gamma_{0,+}\times K_+
  \end{equation*}
  by the de Rham decomposition theorem.
  Hence, we may write~$M=\tilde M/(\Gamma_\pm/\Gamma_{0,\pm})$,
  where~$\tilde M$ is an extra-twisted connected sum as considered above, %
  using only the subgroups~$\Gamma_{0,\pm}\subset\Gamma_\pm$.
\end{rmk}

\subsection{Further questions}

Several questions not answered above are
\begin{itemize}
\item What are $\mu$ and $\xi$ for extra-twisted connected sums?
While $\nu$ was originally defined as a coboundary defect, our computation
for extra-twisted connected sums was analytic. We do not know suitable
coboundaries for extra-twisted connected sums that could be used to
compute $\mu$ and $\xi$.
\item There are now more than~$10^8$ different constructions of
  $G_2$-manifolds. Only a small number of these constructions
  give~$\bar\nu\ne 0$, even fewer give~$3\nmid\bar\nu$.
  In case the number of possible $G_2$-manifolds up to deformation
  is finite, it would be interesting to know if~$\bar\nu=0$
  or~$3|\nu$  %
  is preferred, or if all values occur roughly equally often.
\item Is there a 2-connected 7-manifold that admits a torsion-free
  $G_2$-structure in every
  homotopy class of topological
  $G_2$-structures?
\end{itemize}

\subsubsection*{Acknowledgements}
We thank Jean-Michel Bismut, Uli Bunke, Xianzhe Dai, Matthias Lesch
for inspiring discussions on
adiabatic limits and variational formulas for $\eta$-invariants
on manifolds with boundary.
We thank Alessio Corti, Jesus Martinez Garcia, David Morrison,
Emanuel Scheidegger
and Katrin Wendland for %
helpful information about $K3$-surfaces and Fano threefolds.
We also thank Mark Haskins and Arkadi Schelling for %
talking with us about $G_2$-manifolds and $G_2$-bordism.
We are particularly indebted to Don Zagier for his formula for~$F_{k,\eps}(s)$
in Section~\ref{sect:nuxxtcs}.
SG and JN would like to thank the Simons foundation for its support of their
research under the Simons Collaboration on ``Special Holonomy in Geometry,
Analysis and Physics'' (grants \#488617, Sebastian Goette, and \#488631, Johannes
Nordstr\"om).

\section{Coboundary defects}

Suppose that there is a formula valid for closed $n$-manifolds with
a certain structure, such that each term is well-defined also for $n$-manifolds
with boundary. If each term is additive under gluing %
along %
boundary components,
then the failure of the formula to hold for manifolds
with boundary can be interpreted as an invariant of the boundary itself (with
relevant induced structure).
We explain how combinations of the Hirzebruch signature theorem, the
Atiyah-Singer theorem for the index of the Dirac operator and a relation for
the Euler class of the positive spinor bundle of closed spin 8-manifolds lead
to the definitions of the invariants $\mu$, $\nu$ and $\xi$.

\subsection{Prototypical example}

The first example of such a ``coboundary defect''
invariant is Milnor's $\lambda$-invariant of a closed oriented 7-manifold $M$
with $p_1(M) = 0$.
The starting point in this case is the Hirzebruch signature theorem for a closed
oriented 8-manifold $X$:
\begin{equation}
\label{eq:hirzebruch}
\sigma(X) = \frac{7p_2(X) - p_1(X)^2}{45}
\end{equation}
Here we have implicitly identified $p_2(X), p_1(X)^2 \in H^8(X) \cong \Z$ by
evaluation on the fundamental class.

Now consider instead a compact 8-manifold $W$ with boundary $M$.
Let $H^4_0(W)$ be the image of the push-forward $H^4(W,M) \to H^4(W)$.
For elements $x, y \in H^4_0(W)$ we can define a product $xy \in \Z$ by
picking a pre-image $\bar x \in H^4(W,M)$ of $x$ and setting $xy$ to be
$\bar x y \in H^8(W,M) \cong \Z$.
One makes sense of the signature $\sigma(W)$ as the signature of
this intersection form on $H^4_0(W)$.
If we impose the condition that $p_1(M) = 0$, then $p_1(W) \in H^4_0(W)$,
so $p_1(W)^2 \in \Z$ is well-defined.

According to Novikov additivity~\cite[7.1]{atiyah68}, %
the signature is additive under gluing boundary components:
if $X^8 = W_0 \cup_M W_1$ for manifolds $W_i$ with $\partial W_i = M$ (but
opposite orientations), then $\sigma(X) = \sigma(W_0) + \sigma(W_1)$. %
The integral of %
$p_1^2$ is additive in the same sense.

While there is no good way to interpret
$p_2(W)$ under these conditions, we can eliminate the corresponding term from
\eqref{eq:hirzebruch} by reducing modulo 7:
\begin{equation}
45 \sigma(X) + p_1(X)^2 \equiv 0 \mod 7
\end{equation}
for any closed oriented 8-manifold $X$. The consequence is that if $M$
is a smooth oriented 7-manifold with $p_1(M) = 0$ and $W$ is an oriented
coboundary, then
\[ \lambda(M) := 45 \sigma(W) + p_1(W)^2 \mod 7 \]
in fact depends only on $M$, and not on $W$
(Milnor \cite[Theorem 1]{milnor56}). Because the oriented bordism
group $\Omega^{SO}_7$ is trivial, this allows us to define
$\lambda(M) \in \Z/7$ for any oriented $M$ with $p_1(M) = 0$.

\subsection{The spin characteristic class}
\label{subsec:p}

Since we are interested in invariants of spin manifolds, it will be important
to summarise some properties of the generator of $H^4(B\Spin; \Z)$.
This corresponds to a degree 4 characteristic class $p(E)$ of spin vector
bundles $E$. It is related to the first Pontrjagin class by $p_1(E) = 2p(E)$,
while its mod 2 reduction is the 4th Stiefel-Whitney class $w_4(E)$.

For a manifold $M$, we will abbreviate $p(TM)$ as $p_M$. For a closed
spin manifold of $\dim M = n$, Wu's formula~\cite[Theorem 11.14]{milnor74} implies that $w_4(M)$ coincides %
with the 4th Wu class $v_4(M)$, \ie the Poincar\'e dual to the Steenrod square
$Sq^4 : H^{n-4}(M; \Z/2) \to H^n(M; \Z/2)$.

If $\dim M \leq 7$ then $v_4(M) = 0$, so $p_M$ is even.

If $X$ is closed of dimension 8 then the definition of $v_4(X)$ means it is a
characteristic element for the intersection form on $H^4(X; \Z/2)$, %
that is %
\begin{equation}
\label{eq:p_char}
p_X x = x^2 \mod 2 \quad  \textrm{ for any } x \in H^4(X; \Z) .
\end{equation}
The van der Blij lemma (see Milnor-Husemöller \cite[Chapter II, Lemma 5.2]{milnor73}) implies in turn that
\begin{equation}
\label{eq:mod8}
p_X^2 = \signat(X) \mod 8 .
\end{equation}
One can in fact deduce that \eqref{eq:p_char} and \eqref{eq:mod8} remain valid
also if $X$ is compact with boundary (taking $x \in H^4_0(X;\Z/2)$ in
\eqref{eq:p_char}).
See \cite[\S2.1]{7class} for further details.

\subsection{The Eells-Kuiper invariant and its generalisation}

In the context of closed 7-manifolds that are spin, there is another
relevant formula for closed 8-manifolds in addition to
the signature theorem \eqref{eq:hirzebruch}.
A closed spin $X^8$ has a Dirac operator $D_X$, and by the
Atiyah-Singer theorem its index is computed by the $\ahat$-genus of $X$:
\begin{equation}
\label{eq:index}
\ind D_X = \frac{7p_1(X)^2 - 4p_2(X)}{45 \cdot 2^7}
\end{equation}
While it is possible to define an index of the Dirac operator on a manifold
with boundary, it is not a topological invariant. While we make use of that
below, in the context of defining defect invariants we will need to eliminate
this term. 

To understand how to extract coboundary defect invariants from
\eqref{eq:hirzebruch} and \eqref{eq:index}, it is helpful to rearrange them as
\begin{equation}
\label{eq:pre_ek}
\begin{aligned}
7p_2(X) &= 4p_X^2 + 45\sigma(X) \\
45\cdot2^5\ind D_X + p_2(X) &= 7p_X^2 ;
\end{aligned}
\end{equation}
we have put the terms that have useful interpretations for
manifolds with boundary on the right and the ones that do not on the left,
and we also used $p_1(X) = 2p_X$ to simplify slightly.
Clearly we cannot completely eliminate both $p_2(X)$ and $\ind D_X$ using these
two equations. But if we eliminate the $p_2(X)$ term, then we are left with
\[ 7 \cdot 45 \cdot 2^5 \ind D_X = 45(p_X^2 - \sigma) . \]
Clearly we can eliminate a common factor of 45. However, in view of
\eqref{eq:mod8} it is more natural to reduce to
\[ 28\ind D_X = \frac{p_X^2 - \sigma(X)}{8} . \] 
Thus, if for a closed spin 7-manifold $M$ with $p_M = 0$ we define
\[ \mu(M) := \frac{p_W^2 - \sigma(W)}{8} \in \Z/28 \]
for any spin coboundary $W$, then this will be independent of the choice
of $W$. This is (up to normalisation) the invariant of Eells and
Kuiper \cite{eells62}.
It is the best possible defect invariant that can be extracted from
\eqref{eq:pre_ek} in the following sense:
\begin{itemize}
\item Even given the constraint \eqref{eq:mod8}, $\mu$ can take any
value in $\Z/28$.
\item Given $(p_W^2, \signat(W))$ satisfying \eqref{eq:mod8},
there is a solution $(\ind D_X, p_2(X)) \in \Z^2$ to \eqref{eq:pre_ek}
if and only if $\mu = 0$.
\end{itemize}
If $W$ is a compact spin manifold with boundary $M$ but we drop the condition that $p_M = 0$, then we can no longer interpret $p_W^2$ as a well-defined
element of $\Z$. However, if $p_M$ is divisible by an integer $d$,
then $p_W \mmod d$ belongs to the image of $H^4(W,M;\Z/d)$, so
there is a well-defined $p_W^2 \in H^8(W,M;\Z/d) \cong \Z/d$.
Because $d$ is even, there is also a well-defined Pontrjagin
square in $H^8(W,M;\Z/2d) \cong \Z/2d$. But if we impose that $H^4(M)$
is torsion-free and that there exists $u \in H^4(W)$ such that $du_{|M} = p_M$,
then there is a more elementary way to interpret $p_W^2$ even as an element
of $\Z/2\td$ (for $\td = \lcm(4, d)$ as in the introduction):
if $u'$ is another such element then
\[ (p_W - du')^2 = (p_W - du)^2 + 2dp_W(u'-u) + d^2(u'-u)^2 \in \Z \]
is equal to $(p_W-u)^2$ modulo $2d$ (because $d$ is even) and also modulo
8 (because $p_W$ is a characteristic element for the intersection form %
as explained in %
\eqref{eq:p_char}), so they are equal modulo~$\lcm(8,2d) = 2\td$.

If $H^4(M)$ is torsion-free then one can always find some spin coboundary $W$
and $u \in H^4(W)$ such that $du_{|M} = p_M$. Defining
\begin{equation}
\label{eq:gek_def}
\mu(M) := \frac{(p_W-u)^2 - \sigma(W)}{8} \in \Z/\gcd(28, \tfrac{\td}{4})
\end{equation}
is independent of both $W$ and $u$ (see \cite[Definition 1.8]{7class}).

\subsection{Defect invariants of \texorpdfstring{$G_2$}{G2}-structures}

If we seek invariants of \gtstr s on 7-manifolds rather than just a
spin manifold itself, then one further formula for a closed spin 8-manifold
$X$ becomes relevant.
The integral of the Euler class of the tangent bundle $TX$ is just the Euler
characteristic $\chi(X)$, while the integral of the Euler class of the positive
spinor bundle can be interpreted as the number of zeros (counted with signs)
$n_+(X)$ of any transverse positive spinor field.
They are related by
\begin{equation}
\label{eq:euler}
n_+(X) =  \frac{1}{16}\left(p_1(X)^2 - 4p_2(X) + 8\chi(X) \right) ;
\end{equation}
this appears to have been first established by Gray and Green \cite[p. 89]{gray70}. %
The Euler characteristic of course makes perfect sense also for manifolds
with boundary, and for even-dimensional oriented manifolds 
it is also additive under gluing of boundary components.

On a compact spin manifold $W^8$ with boundary, the number of
zeros of a positive spinor field is not a topological invariant.
However, if we fix a non-vanishing spinor field $s$ on the boundary $M$
and consider transverse positive spinors $\bar s$ on $W$ that restrict to $s$,
then the number of zeros $n_+(W,s)$ does in fact depend only on $s$.
Since a non-vanishing spinor field defines a \gtstr, $n_+(W,s)$ is a sensible
term to consider (only) in the context of manifolds with \gtstr.

We now consider how to define defect invariants from combinations of
\eqref{eq:hirzebruch}, \eqref{eq:index} and \eqref{eq:euler}, which we present
as
\begin{equation}
\label{eq:pre_nuxi}
\begin{aligned}
7p_2(X) &= 4p_X^2 + 45\sigma(X) \\
45\cdot2^5\ind D_X + p_2(X) &= 7p_X^2 \\
p_2(X) &= p_X^2 + 2\chi(X) - 4n_+(X)
\end{aligned}
\end{equation}
In this case we have enough equations that we can eliminate all terms on the
LHS, obtaining
\[ 0 = 7\chi(X) - 14n_+(X) + \frac{3p_X^2-45\signat(X)}{2} \]
Thus for a \gtstr{} on a closed 7-manifold $M$ with $p_M = 0$, defined
by a non-vanishing spinor field $s$, we can define
\[ \xi(s) := 7\chi(W) - 14n_+(W,s) + \frac{3p_W^2 - 45\signat(W)}{2} \in \Z \]
for any spin coboundary $W$.

To capture the remaining constraints from \eqref{eq:pre_nuxi}, there are
many different ways that we could eliminate $p_2(X)$ while leaving an $\ind D_X$
term with a coefficient. If we decide to eliminate the $p_X^2$ term too then
we obtain
\[ -48 \ind D_X = \chi(X) - 2n_+(X) - 3\signat(X) . \]
This allows us to define an invariant of \gstr s by
\begin{equation}
\label{eq:nu_def}
\nu(s) := \chi(W) - 2n_+(W,s) - 3\signat(W) \in \Z/48 .
\end{equation}
Note that $\xi$ and $\nu$ are \emph{not} independent: %
we find using~\eqref{eq:mod8} that %
\begin{equation}
\label{eq:nuxi_rel}
\xi(s) = 7\nu(s) \mod 12
\end{equation}
If we were aiming to identify a ``basic'' set of coboundary defects from
\eqref{eq:pre_nuxi}, we would instead be led to consider $\xi$ together
with a $\Z/4$-valued invariant.

The advantage of instead considering $\nu$ is that it is more robust:
if we drop the condition that $p_M = 0$ then we can no longer define $\xi$
as an integer-valued invariant, but \eqref{eq:nu_def} defines
$\nu(s) \in \Z/48$ for \gtstr s on \emph{any} closed 7-manifold.
However, if we require $H^4(M)$ to be torsion-free then we can define
\[ \xi(s) \in \Z/3\td  \]
analogously to \eqref{eq:gek_def}.

Since we are claiming that $\nu$ and $\xi$ capture all the coboundary-defect
information that can be extracted from \eqref{eq:pre_nuxi}, it should also
be possible to recover $\mu$ from $\nu$ and $\xi$. Indeed it is easy to check
that
\[ 
\frac{\xi(s) - 7\nu(s)}{12} = \mu(M) \mod \gcd \big( 28,\tfrac{\td}{4} \big)\;. \]

\begin{rmk}
In \cite[Definition 1.2 and (10)]{nu}, $\nu$ and $\xi$ are initially defined in
terms of $\Spin(7)$-coboundaries of the \gtstr, \ie using not just a spin
8-manifold $W$ such that $\partial W = M$, but also requiring $W$ to admit a
$\Spin(7)$-structure whose restriction to $M$ is the given \gtstr.
This is equivalent to requiring $n_+(s) = 0$, so leads to a slight
simplification of the defining formulas.
An elementary argument (\cite[Lemma 3.4]{nu}) assures that
$\Spin(7)$-coboundaries exist for any \gtstr.

\label{rmk:g2cob}
One could also ask whether a given \gtstr{} on a 7-manifold admits a
$G_2$-coboundary~$W$. Reducing the structure group of $W$ to $G_2$ defines a
preferred
non-vanishing vector field, so forces $\chi(W) = 0$.
Since $G_2 \subset \Spin(7)$, also $n_+(s) = 0$, so \eqref{eq:nu_def} implies
$\nu = 0 \mmod 3$. In fact, this condition is also sufficient for the existence
of a $G_2$-coboundary \cite{schelling}.
\end{rmk}

\section{Extra-twisted connected sums}\label{Kap3}

We now provide some further details regarding the constructions of twisted
connected sums and extra-twisted connected sums outlined in
\S\ref{Sect:tcs}--\ref{xtcsSect}.

\subsection{ACyl Calabi-Yau manifolds}
\label{subsec:acyl}

The first step in the construction is to produce asymptotically cylindrical
Calabi-Yau 3-folds, \ie complex 3-folds with a complete Ricci-flat Kähler
metric $\omega$ and a choice of (normalised) holomorphic 3-form $\Omega$,
exponentially close to a product structure $(\omega_\infty, \Omega_\infty)$
on an end $\bbrp \times U$.

We will only be concerned with the case when the asymptotic cross-section $U$
is of the form $S^1 \times \kd$, for $\kd$ a K3 surface.
Let $\lnn$ be the circumference of the circle factor, and let $u$ be
a coordinate on $S^1$ with period $\lnn$.
Then there is a \hk triple $(\omega^I, \omega^J, \omega^K)$ on $\kd$ such that
\begin{equation}
\label{eq:cy_limit}
\begin{aligned}
\omega_\infty &= dt \wedge du + \omega^I, \\
\Omega_\infty &= (du - idt) \wedge (\omega^J + i \omega^K)\;.
\end{aligned}
\end{equation}
To produce such ACyl Calabi-Yau manifolds, we use a non-compact version of
Yau's solution of the Calabi conjecture. The following result from \cite{hhn}
is a special case of the Tian-Yau theorem, but with improved control on the
asymptotics.

\begin{thm}
\label{thm:hhn}
Let $Z$ be a closed complex Kähler manifold, and $\kd$ an anticanonical divisor
with trivial normal bundle. Then $Z \setminus \kd$ admits ACyl Ricci-flat
Kähler metrics.
\end{thm}

That $\kd$ is an anticanonical divisor essentially means it is a complex
submanifold Poincar\'e dual to $c_1(Z)$. A convenient way to produce examples
of `building blocks' $Z$ to which Theorem \ref{thm:hhn} can be applied
is to blow up the intersection of two anticanonical divisors in a Fano 3-fold,
\ie a closed complex 3-fold $Y$ where $c_1(Y)$ is a Kähler class.
The topology of such manifolds is well-understood, as is their deformation
theory which is relevant for the matching problem discussed
in~\S\ref{subsec:match}.

\begin{ex}
\label{ex:mm6}
Let $Y \subset \PP^2 \times \PP^2$ be a smooth divisor of bidegree (2,2).
Then the anticanonical bundle $-K_Y$ is the restriction of $\co(1,1)$.
The intersection of two generic anticanonical divisors $\kd_0, \kd_1$
is a smooth curve $C$ of genus 7. Let $Z$ be the blow-up of $Y$ in $C$.
Then the proper transform of $\kd_0$ is an anticanonical divisor in $Z$
with trivial normal bundle.

The `Picard lattice' of $Y$ is $H^2(Y)$ equipped with the bilinear form
$(x,y) \mapsto xy(-K_Y)$. In the basis for $H^2(Y)$ given by the
restrictions of the hyperplane
classes of the $\PP^2$ factors, this form is represented by
\[ \begin{pmatrix} 2 & 4 \\ 4 & 2 \end{pmatrix} . \]
This equals the polarising lattice $N$ of the resulting ACyl Calabi-Yau
3-folds, %
see section~\ref{Sect:tcs}. %
\end{ex}

If $Z$ has a cyclic automorphism group $\Gamma$ that fixes a smooth
anticanonical divisor $\kd$ point-wise, then the restriction of $\Gamma$ to
$Z \setminus \kd$ gives the type of automorphism we need on the ACyl Calabi-Yau.
For instance, in Example \ref{ex:mm6} choosing $Y$, $\kd_0$ and $\kd_1$ to be
invariant under the involution that swaps the $\PP^2$ factors ensures that
this involution lifts to $Z$.

\begin{ex}
\label{ex:cubic}
Let $Y$ be a triple cover of the smooth quadric $Q \subset \PP^4$, branched over
a smooth cubic section $\kd_0 \subset Q$. Let $\kd_1 \subset Y$ be the
pre-image of a generic hyperplane section of $Q$. Then $C := \kd_0 \cap \kd_1$
is a smooth curve of genus 4. Let $Z$ be the blow-up of $Y$ in $C$. Then
the proper transform $\kd \subset Z$ of $\kd_0$ is an anticanonical divisor
with trivial normal bundle, and the branch-switching automorphisms of $Y$
lift to automorphisms of $Z$ that fix $\kd$.

The Picard lattice of $Y$ has rank 1, with a generator that squares to 6.
\end{ex}

\subsection{Gluing ACyl \texorpdfstring{$G_2$}{G2}-manifolds}\relax
\label{subsec:glue}
Choose~$\lnx>0$ and let %
$v$ be a coordinate with period $\lnx$ on~$S^1$.
Given an ACyl Calabi-Yau structure $(\omega, \Omega)$ on $V$, the 3-form
$\varphi := \Re \Omega + dv \wedge \omega$ defines a torsion-free ACyl
\gtstr{} on $S^1 \times V$ and hence a metric with holonomy contained in
(but not equal to) $G_2$; the circumference of the `external' $S^1$ factor
equals $\lnx$.
If the asymptotic limit $(\omega, \Omega)$ is given by \eqref{eq:cy_limit},
then $\varphi$ is asymptotic to
\begin{equation}
\label{eq:limit0}
\phy_\infty = dv \wedge dt \wedge du + dv \wedge \omega^I
+ du \wedge \omega^J + dt \wedge \omega^K
\end{equation}
If $(V,\omega,\Omega)$ admits an isomorphic action by $\Gamma = \Z/k\Z$
with~$k\ge2$ %
as above, then we can extend the action to $S^1 \times V$ by making a 
generator act on the external $S^1$ factor as rotation by
angle~$\frac{2\pi}{k}$;
let $\eps \in \Z/k$ be the unit such that the action of that
generator on the internal~$S^1$ by~$\frac{2\pi\eps}k$.
The ACyl \gtstr{} $\phy$ descends to the quotient
$M := (S^1 \times V)/\Gamma$. It has asymptotic limit $T^2 \times \kd$,
where the $T^2$ factor is isometric to the quotient of a product of
circles of circumference $\lnn$ and $\lnx$ by $\Gamma$.

More precisely, the $T^2$ factor could be described as the quotient of $\C$
by a lattice generated by $\lnx$ and $\frac{\lnx +i \eps\lnn}{k}$, with
a complex coordinate
\begin{equation}
\label{eq:zdef}
z = \anglex + i\anglen \,.
\end{equation}
Given a pair $(\lnn_+, \lnx_+, k_+, \eps_+)$, %
$(\lnn_-, \lnx_-, k_-, \eps_-)$ %
of sets of data defining such tori $T^2_+, T^2_-$ and an angle $\thet \not= 0$,
we consider in the next subsection whether
\[ \C \to \C, \; z \mapsto e^{i\thet}\bar z\] descends to a
well-defined orientation reversing isometry
\begin{equation}
\label{eq:tormat}
\tormat : T^2_+ \to T^2_- \;.
\end{equation}
If we have such a $\thet$, we can attempt to find a diffeomorphism
$\hkr : \kd_+ \to \kd_-$ such that \eqref{eq:cylmap} is an isomorphism
of \gtstr s. Let us now identify this condition in terms of the action
on \hk structures.
In terms of the complex coordinate $z = v+iu$ in \eqref{eq:zdef}, we can
rewrite \eqref{eq:limit0} as
\begin{equation}\label{eq:limit}
 \varphi_\infty = 
\Re \left(dz \wedge (\omega^I - i \omega^J)\right) +
dt \wedge \left(\omega^K - {\textstyle \frac{i}{2}} dz \wedge d\bar z \right).
\end{equation}
For \eqref{eq:cylmap} to be an isomorphism of cylindrical \gtstr s is thus
equivalent to the following.

\begin{dfn}
\label{def:hkr}
Given $\thet \in \R$ and \hk structures
$(\omega^I_\pm, \omega^J_\pm, \omega^K_\pm)$ on K3 surfaces $\kd_\pm$,
call a diffeomorphism $\hkr : \kd_+ \to \kd_-$ a \emph{$\thet$-\hk rotation}
(or simply a \hk rotation if $\thet = \frac{\pi}{2}$) if
\begin{equation}
\label{eq:hkr}
\begin{aligned}
\hkr^*\omega^K_- &= - \omega^K_+ \\
\hkr^*(\omega^I_- + i \omega^J_-) &= e^{i\vartheta} (\omega^I_+ - i\omega^J_+)\;.
\end{aligned}
\end{equation}
\end{dfn}
We will now consider in turn the problems of finding suitable
$\tormat :T^2_+ \to T^2_-$ and $\hkr : \kd_+ \to \kd_-$.

\subsection{Isometries of tori}
\label{subsec:tori}

Given $k_\pm$, identifying the possible data $\eps_\pm$, $\xi_\pm$, $\zeta_\pm$
and $\thet$ for which \eqref{eq:tormat} is well-defined is essentially a
combinatorial problem.
To study it is helpful to associate to such a $\tormat$
a {\em gluing matrix\/}~$G=\bigl(\begin{smallmatrix}m&p\\n&q\end{smallmatrix}\bigr)$
such that
\begin{align}
  \begin{split}\label{eq:mnpqdef}
    \lnx_-\del_{\anglex_-}
    &=\frac1{k_+}\,d\tormat\bigl(m\lnx_+\del_{\anglex_+}+n\lnn_+\del_{\anglen_-}\bigr)\;,\\
    \lnn_-\del_{\anglen_-}
    &=\frac1{k_+}\,d\tormat\bigl(p\lnx_+\del_{\anglex_+}+q\lnn_+\del_{\anglen_-}\bigr)\;.
  \end{split}
\end{align}
Let us write~$s_\pm$ for the ratio~$\frac{\lnx_\pm}{\lnn_\pm}$.
Amongst other relations, the matrix coefficients satisfy
\begin{subequations}\label{eq:mnpqrels}
\begin{gather}
  \det\bigl(\begin{smallmatrix}m&p\\n&q\end{smallmatrix}\bigr)
    =-k_-k_+\;,\label{eq:mnpqdet}\\ %
   mnpq\le 0\;,\label{eq:mnpqprod}\\ %
  \eps_+m-n\equiv\eps_+p-q\equiv 0\mod k_+\;,\label{eq:mnpqeps+}\\
  \eps_-p+m\equiv\eps_-q+n\equiv 0\mod k_-\;,\label{eq:mnpqeps-}
\end{gather}
\end{subequations}
see~\cite{nuxx}.
Because~$\lnx_-\del_{\anglex_-}$ and~$\lnn_-\del_{\anglen_-}$ are perpendicular,
there are three possibilities.
If~$\thet=0$, then~$n=p=0$.
This leads to a manifold with infinite fundamental group,
so we do not consider this case.
If~$\thet=\pm\frac\pi2$, then~$m=q=0$,
and~$\lnx_+=\lnn_-$ and~$\lnn_+=\lnx_-$ are independent of each other.
If~$\thet\notin\frac\pi2\Z$, then~$mnpq<0$, and
\begin{subequations}\label{eq:mnpqrels2}
\begin{gather}
  s_-=\frac{\lnx_-}{\lnn_-}=\sqrt{-\tfrac{mn}{pq}}\;,\qquad
  s_+=\frac{\lnx_+}{\lnn_+}=\sqrt{-\tfrac{nq}{mp}}\;,\label{eq:mnpqs}\\
  \thet=\arg\bigl(ms_++in\bigr)\;.\label{eq:mnpqthet}
\end{gather}
\end{subequations}
For given~$k_\pm$, equations~\eqref{eq:mnpqrels} leave
only finitely many possibilities for~$G$ and~$\eps_\pm$.
If $k_+ = k_- = 1$ then essentially the only
possible gluing matrix is~$G=\bigl(\begin{smallmatrix}0&1\\1&0\end{smallmatrix}\bigr)$, leading to~$\thet = \pm \frac{\pi}{2}$.
\thefigures
\morefigures

For $k_+ = 2$ and $k_- = 1$ there is essentially
only one possibility~$\bigl(\begin{smallmatrix}1&1\\1&-1\end{smallmatrix}\bigr)$.
We can take $\lnn_+ = \lnx_+$ and $\lnn_- = \lnx_-$. That way $T^2_-$ is a
square torus, and $T^2_+$ is a $\Z/2$-quotient of a square torus that is again
a square torus. If we take $\lnn_+ = \sqrt{2} \lnn_-$ then $T^2_+$ and
$T^2_-$ have equal size, and there is an isometry with $\thet = \frac{\pi}{4}$.
To illustrate it we can draw a single lattice corresponding to the two tori
identified by $\tormat$, while adding vectors $\partial_{u_\pm}$ and
$\partial_{v_\pm}$ indicating the directions of the `internal' and `external'
circle directions of the two tori; see Figure~\ref{fig:1/4}.

For $k_+ = k_- = 2$ there are more possibilities, but essentially only two that
lead to simply-connected \gtmfd s. In both of those cases, the tori $T^2_+$ and
$T^2_-$ are `hexagonal'. One possibility is to take
$\lnx_+ = \lnx_- = \sqrt{3}\lnn_+ = \sqrt{3}\lnn_-$, leading to the existence
of an isometry $\tormat$ with $\thet = \frac\pi3$ illustrated in
Figure~\ref{fig:1/3}.
The other has $\lnx_+ = \lnn_- = \sqrt{3}\lnn_+ = \sqrt{3}\lnx_-$ and $\thet = \frac{\pi}{6}$, illustrated in Figure~\ref{fig:1/6}.

Once we allow $k_+$ or $k_-$ to be greater than 2, the number of combinatorial
possibilities increases (while the supply of examples of building blocks with
the relevant symmetry decreases). For $k_+ = 3$ (and $\eps_+ = -1$) and
$k_- = 1$, one possibility is to take
$\lnn_+ = \sqrt{2}\lnx_+ = \sqrt{3}\lnn_- = \sqrt{6}\lnx_-$. That way $T^2_+$
and $T^2_-$ are both rectangular (with %
the proportions of European A4 paper), %
and there is an
isometry $\tormat$ with $\cos \thet = \frac{1}{\sqrt{3}}$ illustrated in
Figure~\ref{fig:A4}; note that $\frac{\thet}{\pi}$ is irrational in this case.

\subsection{The matching problem}
\label{subsec:match}

If we first produce some examples of ACyl Calabi-Yau 3-folds $V_\pm$ with
automorphism groups $\Gamma_\pm$ as in \S\ref{subsec:acyl} and pick a
compatible torus isometry $\tormat$ as in \S\ref{subsec:tori},
it is very unlikely that we will be able to find a $\thet$-\hk rotation between
the asymptotic K3s (for the angle $\thet$ determined by $\tormat$).
A more fruitful approach is to first fix a pair $\calz_+, \calz_-$ of
deformation families of building blocks with automorphism groups $\Gamma_\pm$,
fix $\tormat$, and \emph{then} construct the \emph{pair} $V_+, V_-$ with the
desired $\hkr$ from elements of $\calz_\pm$. 

Fixing a deformation family of blocks $\calz_\pm$ also fixes the polarising
lattice $N_\pm$ of the resulting ACyl Calabi-Yaus. %
Given a $\thet$-\hk rotation $\hkr : \kd_+ \to \kd_-$ between some pair
$(Z_+, \kd_+)$, $(Z_-, \kd_-)$, we can identify both $H^2(\kd_+)$
and $H^2(\kd_-)$ with a fixed copy of the K3 lattice $L$, and hence obtain a
pair of embeddings of $N_+$ and $N_-$ into $L$. %
As in the introduction, section~\ref{Sect:tcs}, %
we refer to this pair as the `configuration' of $\hkr$, and
it controls much of the topology of the resulting \gtmfd s.
It is therefore reasonable to further refine the problem to look for a
$\hkr$ compatible with a fixed configuration.

\begin{rmk}
\label{rmk:nikulin}
According to Nikulin \cite[Theorem 1.12.4]{nikulin79}, an even indefinite
lattice of rank up to 11 has essentially a unique embedding into $L$. As long
as the ranks of $N_+$ and $N_-$ are not too big, specifying a configuration is
therefore essentially equivalent to describing a ``push-out'' %
lattice $W$ that
is spanned by images of isometric embeddings of $N_+$ and~$N_-$.
\end{rmk}

Let us note some necessary conditions on the configuration for the existence of such a $\hkr$.
Observe that $[\omega^I_\pm]$ belongs to
$N^\R_\pm := N_\pm \otimes \R \subset H^2(\kd_\pm;\R)$, and is also the
restriction of a Kähler class from $Z_\pm$. On the other hand, $[\omega^J_\pm]$
and $[\omega^K]_\pm$ are orthogonal to $N^\R_\pm$.
If we let $\pi_\pm : L \to N^\R_\pm$ be the orthogonal projection, then \eqref{eq:hkr} implies that
$\pi_\pm [\omega^I_\mp] = (\cos \thet)[\omega^I_\pm]$, and hence that
$[\omega^I_\pm]$ belongs to the $(\cos \thet)^2$-eigenspace of the self-adjoint
endomorphism $\pi_+\pi_-$ on~$N^\R_\pm$; let us %
denote that by $N^\thet_\pm\subseteq N^\R_\pm$.

Since the positive-definite subspace spanned by $[\omega^I_+]$ and $[\omega^I_-]$ is contained in $W$ while $[\omega^K_+] = -[\omega^K_-]$ is perpendicular
to $W$, %
we see that %
$W$ must be non-degenerate of signature $(2, \rk W - 2)$.

Now let $\Lambda_\pm \subset L$ be the primitive overlattice
of $N_\pm + N_\mp^{\not= \thet}$, where
$N_\mp^{\not=\thet} \subset N_\mp$ is the orthogonal complement
of $N_\mp^\thet$ in $N_\mp$. Recall that the Picard lattice of $\kd_\pm$ is
$H^2(\kd_\pm;\Z) \cap H^{1,1}(\kd_\pm;\R)$. Since $H^{1,1}(\kd_\pm;\R)$ is the
orthogonal complement in $H^2(\kd_\pm;\R)$ to the span of
$[\omega^J_\pm]$ and $[\omega^K_\pm]$, \eqref{eq:hkr} further forces that
$\Lambda_\pm$ is contained in the Picard lattice
of $\kd_\pm$, so ``$\kd_\pm$~is $\Lambda_\pm$-polarised''. %

In summary, given a pair of primitive embeddings $N_+, N_- \into L$ of the polarising lattices of a pair of deformation families $\calz_+, \calz_-$ of building blocks, two necessary conditions for finding a $\thet$-\hk rotation between
asymptotic K3s in some ACyl Calabi-Yau 3-folds arising from some elements of
$\calz_+$ and $\calz_-$ are that:
\begin{enumerate}
\item $W := N_+ + N_-$ is non-degenerate of signature $(2,\rk W -2)$
\item $N^\thet_\pm$ contains the restriction of some Kähler class from $Z_\pm$;
in particular $N^\thet_\pm$ is non-trivial
\item there are some elements $(Z_\pm, \kd_\pm) \in \calz_\pm$ such that
$\kd_\pm$ is $\Lambda_\pm$-polarised.
\end{enumerate}
On the other hand, a combination of the Torelli theorem with a
more precise statement of Theorem \ref{thm:hhn} turns out to show that a
sufficient condition for finding a $\thet$-\hk rotation compatible with
the configuration is given essentially by (i) and (ii) together with
\begin{enumerate}
\item[(iii')] a \emph{generic} element of the moduli space of $\Lambda_\pm$-polarised
K3s appears as the anticanonical divisor in some element of $\calz_\pm$
\end{enumerate}
A general principle is that a generic $N_\pm$-polarised K3 surface does appear
as an anticanonical divisor in some element of $\calz_\pm$.
For example, for blocks obtained from Fano 3-folds, %
as in %
Example \ref{ex:mm6}, this is a consequence of the results of Beauville \cite{beauville04} on the deformation theory of anticanonical divisors in Fano 3-folds
(see \cite[Proposition 6.9]{cym}).
The matching problem is therefore easiest to solve if one restricts attention
to configurations where $\Lambda_\pm = N_\pm$.
That is equivalent to requiring
that the only configuration angles in~\eqref{eq:mrho} are $0$ and $\pm 2\thet$.

\subsection{Examples of matchings}

For $\thet = \frac{\pi}{2}$, it is very easy to produce such configurations
where $\Lambda_\pm = N_\pm$:
simply take the push-out $W$ of Remark \ref{rmk:nikulin} to be the
perpendicular direct sum $N_+ \perp N_-$; then (i) and (ii) are automatically
satisfied too. This way one can produce literally millions of matchings,
see \cite{g2m}. However, there is limited diversity among the topological types
realised this way, \eg they all have $\mu = 0$ \cite[Corollary 3.7]{exotic}.

On the other hand, if $\thet \not= \frac{\pi}{2}$, then for a given pair of
polarising lattices $N_+$ and $N_-$ there need not be any
configurations at all with $N_\pm = \Lambda_\pm$.
For polarising lattices of rank 1, it is not so difficult to decide whether
such a configuration exists.

\begin{ex}
\label{ex:irrat}
Let $\calz_+$ be the deformation family of blocks with automorphism group
$\Gamma_+ \cong \Z/3$ described in Example \ref{ex:cubic}, and let
$\calz_-$ be the family of blocks obtained from blow-ups of Fano 3-folds
of rank 1, index 1 and degree 2 (\cf \cite[Example $7.1^1_2$]{g2m} %
in the notation used there). %
The relevant polarising lattices are $N_+ = (6)$ and $N_- = (2)$.
By the reasoning in Remark \ref{rmk:nikulin}, the matrix
\[ W = \begin{pmatrix} 6 & 2 \\ 2 & 2 \end{pmatrix} \]
defines a configuration of $N_+$ and $N_-$. The angle $\thet$ between the
basis vectors has
\[ (\cos \thet)^2 = \frac{2^2}{2\cdot 6} = \frac{1}{3}\;. \]
We can find a $\thet$-\hk rotation compatible with this configuration, and
hence form an extra-twisted connected sum using the torus matching illustrated
in Figure~\ref{fig:A4}.
\end{ex}

On the other hand, for polarising lattices of higher rank the existence
can be less immediately obvious.

\begin{ex}
\label{ex:rk2pi4}
Let $\calz_+$ be the deformation family of blocks with automorphism group
$\Gamma_+ \cong \Z/2$ described in Example \ref{ex:mm6}, and let
$\calz_-$ be the family of blocks obtained by blowing up the blow-up of $\PP^3$ in a conic (number 30 in the Mori-Mukai classification of rank 2 Fano 3-folds,
see \cite[Entry 30 of Table 3]{exotic}).
The relevant polarising lattices are
\[ N_+ = \begin{pmatrix} 2 & 4 \\ 4 & 2 \end{pmatrix} , \qquad
N_- = \begin{pmatrix} 6 & 6 \\ 6 & 4 \end{pmatrix} . \]
Then
\[ W = \begin{pmatrix} 2 & 4 & 3 & 4 \\ 4 & 2 & 3 & 2 \\
3 & 3 & 6 & 6 \\ 4 & 2 & 6 & 4 \end{pmatrix} \]
defines a configuration of $N_+$ and $N_-$,
such that $N_\pm^{\frac{\pi}{4}} = N_\pm$.
We can find a $\frac{\pi}{4}$-\hk rotation compatible with this configuration,
and hence form an extra-twisted connected sum using the torus matching
illustrated in Figure~\ref{fig:1/4}.
\end{ex}

Even if we look for configurations without the assumption that
$N_\pm = \Lambda_\pm$, the conditions (i) and (ii) on their own can be still
be quite restrictive. But having found such a configuration, one then has to
check condition (iii'). This typically requires some detailed understanding
of the particular families of building blocks involved.

\begin{ex}
\label{ex:twistedcubic}
Take both $\calz_+$ and $\calz_-$ to be the family of blocks obtained from
blowing up the blow-up of $\PP^3$ in a twisted cubic (number 27 in the
Mori-Mukai classification of rank 2 Fano 3-folds,
see \cite[Entry 27 in Table 3]{exotic}).
In this case the polarising lattices are
\[ N_+ = N_- = \begin{pmatrix} 4 & 5 \\ 5 & 2 \end{pmatrix} . \]
We can define a configuration satisfying condition (i)
(with $\thet = \frac{\pi}{2}$) and (ii)
using the push-out
\[ W =  \begin{pmatrix} 4 & 5 & 1 & -1 \\ 5 & 2 & -1 & 1
\\ 1 & -1 & 4 & 5 \\ -1 & 1 & 5 & 2 \end{pmatrix} . \]
Now $\Lambda_\pm$ is a rank 3 overlattice of $N_\pm$, with quadratic form
represented by
\[ \begin{pmatrix} 4 & 5 & 16 \\ 5 & 2 & -16 \\ 16 & -16 & -272 \end{pmatrix}
. \]
It is checked in \cite[Lemma 7.7]{exotic} that any K3 surface with Picard
lattice isomorphic to that can be embedded as an anticanonical divisor in the
blow-up of $\PP^3$ in a twisted cubic, so that (iii') holds.
Thus it is possible to form a rectangular twisted connected sum of
two blocks from this family. The resulting \gtmfd s have $b_3 = 101$, $d = 8$
and $\mu = 1$, and are used in Example~\ref{ex:exotic}.
\end{ex}

\section{The extended \texorpdfstring{$\nu$}{nu}-invariant}\label{NuSect}
By definition, coboundary defect invariants for $M$ 
can be computed if one knows enough about some appropriate manifold~$W$ with~$\del W=M$.
For rectangular twisted connected sums, this was used in~\cite{nu}
to show that~$\nu(s)=24$, in \cite{exotic} to compute the generalised
Eells-Kuiper invariant, and recently by Wallis to
compute $\xi(s)$ \cite{wallis18}.
For extra-twisted connected sums, zero-bordisms are harder to find, and %
we therefore pursue a different approach to computing $\nu$.

We rewrite the definition of~$\nu(M)$ using the Atiyah-Patodi-Singer index
theorem for manifolds with boundary and Mathai-Quillen currents.
This yields a formula for~$\nu$ in terms of $\eta$-invariants and Mathai-Quillen
currents.
In the case of $G_2$-holonomy, the Mathai-Quillen terms drop out,
and the $\eta$-invariants become $\R$-valued rather than just~$\R/2\Z$-valued.
This way, the $\nu$-invariant lifts to a $\Z$-valued invariant~$\bar\nu$,
the {\em extended $\nu$-invariant,\/}
that is locally constant on the moduli space of $G_2$-holonomy manifolds.
It is possible to compute $\bar\nu$ for extra-twisted connected sums, %
see examples~\ref{ex36} and~\ref{ex48} above.

\subsection{The analytic description of the \texorpdfstring{$\nu$}{nu}-invariant}\label{APSsect}
The definition of~$\nu(s)$ in~\eqref{eq:nu_def} involves
the signature~$\sigma(X)$ of an $8$-manifold~$X$,
which can be written as the analytic index
of the signature operator~$B_X$ on~$X$.
Implicitly, $\nu(s)$ also involves the index
of the Atiyah-Singer spin Dirac operator~$D_X$ on~$X$.

The Atiyah-Patodi-Singer index theorem allows us to write~$\signat(W)$
as an analytic index of the signature operator~$B_W$ on an $8$-manifold
with boundary~$\del W=M$.
We assume that~$W$ has product geometry near its boundary.
Let~$\nabla^{TW}$ be the Levi-Civita connection,
and let~$L\bigl(TW,\nabla^{TW}\bigr)\in\Omega^\bullet(W)$
be the Chern-Weil representative of the $L$-class.
If~$B_M$ denotes the odd signature operator on the boundary~$M$,
with spectrum~$\cdots\le\lambda_0\le\lambda_1\le\cdots$
counted with multiplicities,
we can define its {\em $\eta$-invariant\/} by
\begin{equation*}
  \eta(B_M)=\sum_{\lambda_i\ne 0}\sign(\lambda_i)\,\abs{\lambda_i}^{-s}\Bigr|_{s=0}
  =\int_0^\infty\tr\bigl(B_M\,e^{-tB_M^2}\bigr)\,\frac{dt}{\sqrt{\pi t}}\;.
\end{equation*}
The spectral expression is defined if the real part of~$s$ is sufficiently
large and has a meromorphic continuation %
that is holomorphic at~$s=0$.
The $\eta$-invariant is the value at~$s=0$,
or equivalently the value of the integral on the right.
The Atiyah-Patodi-Singer signature theorem~\cite[Thm~4.14]{atiyah75}
implies that
\begin{equation}\label{eq:APS_sign}
  \sigma(W)
  =\int_WL\bigl(TW,\nabla^{TW}\bigr)-\eta(B_M)\;.
\end{equation}

Similarly,
let~$D_W$ denote the spin Dirac operator on~$W$ with the given spin structure,
and let~$D_M$ denote the spin Dirac operator on~$M$.
Let~$\ind{\mathrm{APS}}(D_W)\in\Z$ denote the analytic index
of the spin Dirac operator with respect to the
Atiyah-Patodi-Singer boundary conditions,
let~$\hat A(TW,\nabla^{TW})$ be the Chern-Weil representative
of the $\hat A$-class, let~$\eta(D_M)$ be the defined as above,
and let~$h(D_M)=\dim\ker(D_M)$.
The Atiyah-Patodi-Singer index theorem~\cite[Thm~4.2]{atiyah75} states that
\begin{equation}\label{eq:APS_dirac}
  \ind_{\mathrm{APS}}(D_W)=\int_W\hat A\bigl(TW,\nabla^{TW}\bigr)
  -\frac{\eta+h}2(D_M)\;.
\end{equation}

The Euler class of the positive spinor bundle has to be treated differently.
Let~$\pi\colon E\to W$ be a Euclidean vector bundle
with metric~$g^E$ and compatible connection~$\nabla^E$.
Mathai and Quillen \cite{mathai86} defined a current~$\psi(\nabla^E,g^E)$
on the total space~$E$,
which is singular along the zero section~$W\subset TW$,
such that
\begin{equation*}
  d\psi(\nabla^E,g^E)
  =\pi^*e\bigl(E,\nabla^E\bigr)-\delta_W\;.
\end{equation*}
Here, $e(E,\nabla^E)$ is the Euler class of~$E$
and~$\delta_W$ denotes the Dirac delta distribution on~$TW$
along the zero section~$W$.
As a bundle~$E$, we consider the positive spinor bundle~$S^+W\to W$,
so~$SM=S^+W|_M$ is the spinor bundle on~$M$.
If~$\bar s\in\Gamma(S^+W)$ extends a nowhere vanishing spinor~$s$ on~$M$,
then~\cite[Thm~7.6]{mathai86}, see also~\cite[Thm~3.7]{BZtor}, implies
\begin{equation}\label{eq:MQ}
  n_+(W,s)
  =\int_W\bar s^*\delta_W
  =\int_We\bigl(S^+W,\nabla^{S^+W}\bigr)
  -\int_Ms^*\psi\bigl(\nabla^{SM},g^{SM}\bigr)
\end{equation}
by Stokes' theorem.
Thus, at least formally, the integral of the Mathai-Quillen form over~$M$
is analogous to the $\eta$-invariants in~\eqref{eq:APS_sign}
and~\eqref{eq:APS_dirac}.
We combine~\eqref{eq:euler} and~\eqref{eq:nu_def}
with~\eqref{eq:APS_sign}--\eqref{eq:MQ}
to get an intrinsic formula for the $\nu$-invariant.

\begin{thm}\label{nu_thm}
  Let~$s\in\Gamma(SM)$ define a $G_2$-structure on a spin 7-manifold~$M$. Then
  \begin{equation*}
    \nu(s)=3\,\eta(B_M)
	-24\,(\eta+h)(D_M)+2\int_Ms^*\psi\bigl(\nabla^{SM},g^{SM}\bigr)\in\Z/48\;.\qed
  \end{equation*}
\end{thm}

\subsection{The extended \texorpdfstring{$\nu$}{nu}-invariant}\label{nubarSect}
Let us now assume that~$(M,g)$ has holonomy~$G_2$.
Then the defining spinor~$s\in\Gamma(SM)$ is parallel,
and~$s^*\psi(g^{SM},\nabla^{SM})$ vanishes by construction,
\mbox{see~\cite[Lemma~1.3]{eta}}.
Hence, Theorem~\ref{nu_thm} becomes
\begin{equation}\label{eq:nubar}
  \nu(s)=3\,\eta(B_M)
		-24\,(\eta+h)(D_M)\in\Z/48\;.
\end{equation}

We recall that the $\eta$-invariants~$\eta(B_M)$, $\eta(D_M)$ depend
on the spectrum of~$B_M$ and~$D_M$, and hence on the Riemannian geometry
of~$(M,g)$. If one varies the metric~$g$, the corresponding variation
formula for $\eta$-invariants typically contains two terms.
The first term is an integral of a Chern-Simons class over~$M$,
which varies continuously in~$g$.
Since~$\nu(s)$ is always an integer,
the variation terms for the two $\eta$-invariants involved must cancel
for families of metrics with holonomy in~$G_2$.

The second term is a $\Z$-valued spectral flow,
so the $\eta$-invariant, %
or more precisely the expression~$\frac{\eta+h}2$, %
can jump by integers.
However, spectral flow can only occur if eigenvalues
of the relevant operator change sign.
In this case, the dimension~$h$ of the kernel must change.
The kernel of~$B_M$ describes de Rham cohomology,
so~$h(B_M)$ is constant and~$\eta(B_M)$ never jumps.
For the spin Dirac operator this is false in general;
this gives an alternative explanation why~$\nu(s)$ takes values in~$\Z/48$
and not in~$\Z$.

However, if the holonomy group of~$(M,g)$ is a subgroup of~$G_2$,
then~$(M,g)$ is Ricci flat.
The Lichnerowicz formula becomes~$D_M^2= %
(\nabla^{SM})^* %
\nabla^{SM}$.
Because~$M$ is closed, this implies that every harmonic spinor is parallel.
If the holonomy group of~$M$ is the full group~$G_2$,
then the space of parallel spinors is spanned by the defining spinor~$s$,
so we have~$h(D_M)=%
1$.
Otherwise, by Ricci flatness, the entire first de Rham cohomology
can be represented by parallel 1-forms,
and Clifford multiplication~$c_{\,\cdot\,}s$ gives an isomorphism
from~$H^1(M;\R)$ to the subspace of parallel spinors perpendicular to~$s$.
Hence~$h(D_M)=1+b_1(M)$ is constant on the moduli space
of $G_2$-holonomy metrics, and the spin Dirac operator has no spectral flow.
Therefore, the right hand side of~\eqref{eq:nubar} is locally constant
on the $G_2$-moduli space.

\begin{dfn}[{\cite[Definition~1.4]{eta}}]\label{def:nubar}
  For a closed Riemannian 7-manifold~$(M,g)$ with holonomy
  contained in~$G_2$, put
  \begin{equation*}
    \bar\nu(M,g)=3\eta(B_M)-24\eta(D_M)\;.
  \end{equation*}
\end{dfn}

\begin{cor}\label{cor:nubar}
  For a closed Riemannian 7-manifold~$(M,g)$ with holonomy
  contained in~$G_2$ and with defining parallel spinor~$s$, we have
  \begin{equation*}
    \nu(s)=\bar\nu(M,g)-24\,(1+b_1(M))\mod 48\;. %
  \end{equation*}
\end{cor}

One could argue that we should have changed either~\eqref{eq:nu_def}
or Definition~\ref{def:nubar} in order to avoid the correction term
$24(1+b_1(M))$.
But both definitions are the most natural in their respective realm.
In particular, $\bar\nu(M,g)$ changes sign under reversing the orientation
of $M$, and so vanishes if~$(M,g)$ admits an orientation
reversing isometry.

\subsection{Extra-twisted connected sums}\label{sect:nuxxtcs} %
We return to extra-twisted connected sums and sketch a proof %
of Theorem~\ref{thm:involutions}.
Let~$M_\ell=M_{+,\ell}\cup M_{-,\ell}$ as in Section~\ref{Kap3}
be such that~$Y=M_{+,\ell}\cap M_{-,\ell}=\del M_{+,\ell}=\del M_{-,\ell}$,
and
\begin{equation*}
  M_{\pm,\ell}
  =\bigl(S^1\times(V_\pm\setminus((\ell,\infty)\times S^1\times\Sigma_\pm))\bigr)
  \bigm/\Gamma_\pm
\end{equation*}
with~$\Gamma_\pm\cong\Z/k_\pm$.
The parameter~$\ell$ stands for the length of the cylindrical neck.
There is a closed \gtstr~$\phy_\ell$ on~$M_\ell$,
and a torsion free \gtstr~$\bar\phy_\ell$ nearby.

We apply the $\R$-valued gluing formula for $\eta$-invariants
by Bunke~\cite{bunke95} and Kirk-Lesch~\cite{KiLe}.
In~\cite{eta}, we construct operators~$D_{M,\ell}$
and~$B_{M,\ell}$ that are of product type on a neighbourhood of~$Y$,
and have the same kernels as the corresponding operators
on the $G_2$-manifold~$(M_\ell,\bar\phy_\ell)$.
The harmonic spinors on~$Y$ that extend to harmonic spinors
of the restrictions~$D_{M_\pm,\ell}$, $B_{M_\pm,\ell}$
to~$M_{\pm,\ell}$ form Lagrangian subspaces~$L_{D_\pm}\subset\ker(D_Y)$
independent of~$\ell$.
Similarly, harmonic forms
representing~$\im\bigl(H^\bullet(M_\pm;\R)\to H^\bullet(Y;\R)\bigr)$
form Lagrangians~$L_{B_\pm}\subset\ker(B_Y)$.
We modify the APS boundary conditions for the operators~$D_{M_\pm}$
and~$B_{M_\pm}$ on the two halves~$M_\pm$ by these Lagrangian subspaces
and define~$\eta_{\mathrm{APS}}(D_{M_\pm};L_{D_\pm})$
and~$\eta_{\mathrm{APS}}(D_{M_\pm};L_{B_\pm})$
with respect to those boundary conditions.

Recall the polarising lattices~$N_\pm$ inside the K3 lattice~$L$
from Section~\ref{Sect:tcs}.
Let~$A_\pm$ denote the reflections of~$L\otimes\R=H^2(\Sigma;\R)$
in the subspaces~$N_\pm$.
Then the \emph{configuration angles} are the arguments~$\alpha^+_1,
\alpha^+_2, \alpha^+_3$ and
$\alpha^-_1, \ldots, \alpha^-_{19}$ of the eigenvalues of the restrictions of
$A_+ \circ A_-$ to an invariant positive or negative subspace
of~$H^2(\Sigma;\R)$, respectively.
We always have~$\{\alpha^+_1,\alpha^+_2,\alpha^+_3\}=\{0,\pm2\thet\}$.
We define
  \begin{equation}
  \label{eq:mrho}
    m_\rho(L;N_+,N_-)
    =\Sign\rho\,\Bigl(\#\bigl\{\,j
		\bigm|\alpha_j^-\in\{\pi-\abs\rho,\pi\}\,\bigr\}-1
        +2\,\#\bigl\{\,j
		\bigm|\alpha_j^-\in(\pi-\abs\rho,\pi)\,\bigr\}\Bigr)\;.
  \end{equation}
  By~\cite{bunke95} and~\cite{KiLe}, see~\cite[Thm~1]{eta},
  we find that
\begin{align*}
  \bar\nu(M)
  &=\bar\nu(M_+)+\bar\nu(M_-)-72\frac\rho\pi+3m_\rho(L;N_+,N_-)\;,\\
  \text{where}\qquad
  \bar\nu(M_\pm)&=\lim_{\ell\to\infty}\bigl(3\eta_{\mathrm{APS}}(B_{M_\pm,\ell};L_{B_\pm})
    -24\eta(D_{M_\pm,\ell};L_{D_\pm})\bigr)\;.
\end{align*}

To describe the remaining ingredients of Theorem~\ref{thm:involutions},
let~$\lnn_\pm$ and~$\lnx_\pm$ denote the lengths of the ``interior''
and ``exterior'' circle factors as in Section~\ref{Kap3},
and define~$s_\pm$ as in~\eqref{eq:mnpqs}.
We will now set the exterior radius to~$\lnx_\pm=a\lnn_\pm$ instead
and consider~$M_{\pm,a}=(S^1_{a\lnn_\pm}\times V_\pm)/\Gamma_\pm$.
To compute~$\bar\nu(M_{\pm,a})$,
we will compute its limit as~$a\to 0$,
and the variation of~$\bar\nu(M_{\pm,a})$ as~$a$ changes.

To describe the limit~$a\to 0$,
let~$\gamma_\pm\in\Gamma_\pm$ be the generator that rotates the exterior
circle factor by~$\frac{2\pi}{k_\pm}$.
Let~$V_\pm^{0,j}\subset V_\pm$ be the set of isolated %
fixed points %
of~$\gamma_\pm^j$, and for~$p\in V_\pm^{0,j}$,
let~$\alpha_{j,1}(p)$, $\alpha_{j,2}(p)$, $\alpha_{j,3}(p)$
denote the angles of the $\gamma_\pm^j$-action on~$T_pV_\pm$.
Because the $\Gamma_\pm$-action preserves the holomorphic volume form,
these angles can be chosen such that their sum is~$0$.
Then the isolated %
fixed points %
contribute to~$\bar\nu(M_\pm)$ by
\begin{equation}
\label{eq:dedekind}
  D_{\gamma_\pm}(V_\pm)
  =\lim_{a\to0}\bar\nu(M_{\pm,a})
  =\frac3{k_\pm}\sum_{j=1}^{k_\pm-1}\cot\frac{\pi j}{k_\pm}
  \sum_{p\in V^{0,j}_\pm}
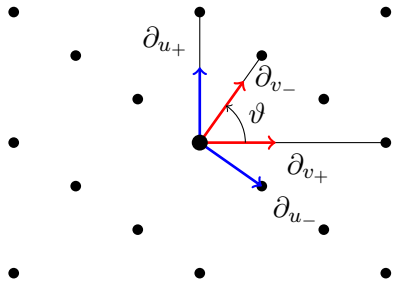
  \frac{\cos\frac{\alpha_{j,1}(p)}2\cos\frac{\alpha_{j,2}(p)}2\cos\frac{\alpha_{j,3}(p)}2-1}
       {\sin\frac{\alpha_{j,1}(p)}2\sin\frac{\alpha_{j,2}(p)}2\sin\frac{\alpha_{j,3}(p)}2}\;,
\end{equation}
see~\cite%
{nuxx}.
This is proved using methods from~\cite{goette11}.

Another contribution arises as a boundary term in the variational formula
for $\eta$-invariants on manifolds with boundary
by Bismut-Cheeger~\cite{BCh4} and Dai-Freed~\cite{daifreed}.
Assume that the generator~$\gamma_\pm$ of~$\Gamma_\pm$
rotates the interior circle by an angle~$\frac{2\pi\eps_\pm}{k_\pm}$
as above.
Let~$\sigma_{-1}(n)=\sum_{d\mid n}d^{-1}$,
and let~$L(\tau)$ denote the logarithm of the Dedekind $\eta$-function,
defined for~$\tau\in\mathcal H\subset\C$ in the upper half plane by
\begin{equation*}
  L(\tau)=\frac{\pi i\tau}{12}-\sum_{n=1}^\infty\sigma_{-1}(n)\,e^{2\pi in\tau}\;.
\end{equation*}
Then the last contribution to~$\bar\nu(M)$ is
\begin{equation}
\label{eq:F}
\begin{aligned}
  F_\pm&=\int_0^{s_\pm}\frac d{da}\bar\nu(M_{\pm,a})
  =\frac{144}{\pi}\,F_{k_\pm,\eps_\pm}(s_\pm)\;,\\
  \text{where}\qquad
  F_{k,\eps}(s)
  &=
  iL\biggl(\frac{si+\eps}{k}\biggr)
  -iL\biggl(\frac{si-\eps}{k}\biggr)
  +c_{k,\eps}\;,
\end{aligned}
\end{equation}
see~\cite%
{nuxx}.
The constant~$c_{k,\eps}$ takes the special values
\begin{equation}\label{eq:ckeps}
    c_{k,\eps}=
    \begin{cases}
      -\eps\pi\,\frac{k^2-3k+1}{6k}	&\text{if~$\eps=\pm 1$, and}\\
      \frac{\pi\eps}{6k}		&\text{if~$\eps^2\equiv-1$ modulo~$k$.}
    \end{cases}
\end{equation}
We are grateful to Don Zagier for the formulas above
for~$F_{k,\eps}(s)$ and~$c_{k,\eps}$.

The explicit values of~$L$ are hard to determine.
Instead, one may use the functional equations
\begin{equation}\label{eq:fneq}
  L(\tau+1)=\frac{\pi i}{12}+L(\tau)\qquad\text{and}\qquad
  L\biggl(-\frac1\tau\biggr)=\frac12\,\log\Bigl(\frac\tau i\Bigr)+L(\tau)
\end{equation}
to compute the sum of all values of~$L$
occurring in Theorem~\ref{thm:involutions} for a particular
extra-twisted connected sum. %

\begin{ex}\label{ex:xbord}
  We consider Example~\ref{ex:irrat}, where~$k_+=3$, $k_-=1$.
  By construction in Example~\ref{ex:cubic}, the group~$\Gamma_+$
  acts without isolated %
  fixed points %
  on~$V_+$, so we have~$D_{\gamma_+}(V_+)=0$.
  And because~$k_-=1$, also~$D_{\gamma_-}(V_-)=0$.
  
  From the gluing matrix~$G=\gmatrix mpnq=\gmatrix112{-1}$
  in Figure~\ref{fig:A4} we conclude
  that~$\eps_+=-1$, $s_+=\sqrt 2=s_-$.
  Because~$k_-=1$, we have~$F_{k_-,\eps_-}(s_-)=0$.
  Using~\eqref{eq:ckeps} and~\eqref{eq:fneq}, we compute
  \begin{align*}
    F_{k_+,\eps_+}(s_+)
    &=iL\biggl(\frac{\sqrt 2i-1}3\biggr)
    -iL\biggl(\frac{\sqrt 2i+1}3\biggr)+c_{3,-1}\\
    &=\frac i2\,\log\frac{\sqrt 2-i}{\sqrt 2+i}
    +iL\bigl(\sqrt 2i+1\bigr)-iL\bigl(\sqrt 2i-1\bigr)+\frac\pi{18}\\
    &=\frac i2\,\log\frac{1-\sqrt 8i}3
    -\frac\pi6+\frac\pi{18}
    =\frac12\arc\cos\frac13-\frac\pi9\;.
  \end{align*}

  Because both~$N_+$ and~$N_-$ have rank~$1$,
  both lie in~$H^{2,+}(\Sigma;\R)$.
  So~$A_+\circ A_-$ acts as the identity on~$H^{2,-}(\Sigma;\R)$,
  and hence~$\alpha^-_1=\cdots=\alpha^-_{19}=0$.
  The angle~$\thet=\arc\cos\frac1{\sqrt3}$ is acute, so~$\rho>0$, %
  hence~$m_\rho(L;N_+,N_-)=-1$.
  Combining all this information, Theorem~\ref{thm:involutions} gives
  \begin{equation*}
    \bar\nu(M)=\frac{144}\pi\biggl(\frac12\arc\cos\frac13-\frac\pi9\biggr)
    -\frac{72}\pi\,\biggl(\pi-2\arc\cos\frac1{\sqrt 3}\biggr)-3=-19\;.
  \end{equation*}
  We see that~$3\nmid\bar\nu(M)$,
  so~$(M,g)$ is indeed not $G_2$-nullbordant.
\end{ex}

\subsection{Elementary hyperbolic geometry}\label{sec:hyp}
There is an alternative way to treat %
the %
variational term~$F_++F_-$. We can compute it as the area of a certain
ideal hyperbolic polygon, see~\cite%
{nuxx}.
To this end, we regard the upper half plane~$\mathcal H$ as space
of conformal structures on a fixed torus.
Then~$\mathcal H$ carries a tautological family of flat tori.
Let~$\tilde\eta(\mathbb A)\in\Omega^1(\mathcal H)$ be the $\eta$-form
of the spin Dirac operator of this family.
Using the variation formula for $\eta$-invariants
on manifolds with boundary in~\cite{BCh4} and~\cite{daifreed},
we represent~$F_\pm$ as
\begin{equation}\label{eq:fpmeta}
  F_\pm=\pm 288\int_{\gamma_\pm}\tilde\eta(\mathbb A)\;.
\end{equation}
Using local index theory, one expresses the exterior derivative
of the $\eta$-form in terms of the hyperbolic volume form~$dA_{\mathrm{hyp}}$ as
\begin{equation}\label{eq:deta}
  d\tilde\eta(\mathbb A)=\frac1{4\pi}\,dA_{\mathrm{hyp}}\;.
\end{equation}

Let~$\gamma_\pm\colon(0,s_\pm]\to\mathcal H$ 
represent the families~$(S^1_{a\lnn_\pm}\times S^1_{\lnn_\pm})/\Gamma_\pm$.
Then~$\gamma_\pm$ are hyperbolic rays.
As we explain in~\cite%
{nuxx}, the ray~$\gamma_+$ goes
from~$\frac{\eps_+}{k_+}\in\R\cup\{\infty\}=\del_\infty\mathcal H$
vertically to the point~$\frac{\eps_++is_+}{k_+}$ representing~$T^2_+$.
The ray~$\gamma_-$ goes from~$\frac{\eps_+}{k_+}-\frac n{k_+m}$
to~$\frac{\eps_++is_+}{k_+}$ along a hyperbolic geodesic
with second endpoint~$\frac{\eps_+}{k_+}-\frac q{k_+p}$.
We can now complete~$\gamma_+\cup\gamma_-$ to an ideal hyperbolic polygon~$P$
of finite area using geodesics along which~$\tilde\eta(\mathbb A)$
vanishes for symmetry reasons;
these are hyperbolic geodesics joining points~$\frac ab$,
$\frac cd\in\Q$ with~$k=\abs{ad-bc}\in\{1,2\}$,
corresponding to families of rectangular~($k=1$) and rhombic~($k=2$) tori,
respectively.

By Stokes theorem and~\eqref{eq:fpmeta} and~\eqref{eq:deta},
we can express~$F_++F_-$ as the sum of~$\frac{72}\pi\,A_{\mathrm{hyp}}(P)$
and contributions from the cusps of~$P$.
Using a strict version of the adiabatic limit formula for families by Bunke and
Ma~\cite{BuMa}, a cusp at~$\frac ef$ between geodesics to~$x$
and~$y\in\del_\infty\mathcal H$ contributes to~$F_++F_-$
by~$-24\measuredangle_{\frac ef}(x,y)$, where the {\em cusp angle\/} is given as
\begin{equation}\label{eq:cuspangle}
  \measuredangle_{\frac ef}(x,y)=\frac{x-y}{(fx-e)(fy-e)}\in\R
\end{equation}
if~$\frac ef$ is a reduced fraction.
Recall that the hyperbolic area of a polygon can be computed
from its angles and the number of corners.
Because the rays~$\gamma_\pm$ meet at angle~$2\thet$,
this approach explains in particular why the final value of~$\bar\nu(M)$
is rational even though the terms~$-72\frac\rho\pi$
and~$F_\pm$
in Theorem~\ref{thm:involutions} can be irrational for~$k_+>2$ or~$k_->2$.

\hypfigure

\begin{ex}\label{ex:hyp}
  We still consider the example above, but compute~$F_++F_-$ using hyperbolic
  geometry.
  Here, $\gamma_+$ lies on the vertical line
  with real part~$-\frac13$,
  and the ray~$\gamma_-$ lies on the hyperbolic geodesic from~$-1$ to~$0$.
  We complete to a hyperbolic polygon with another cusp at~$-\frac12$,
  see Figure~\ref{fig:hyp}.
  Because~$P$ consists of two ideal triangles,
  we have~$A_{\mathrm{hyp}}(P)=2\pi-2\thet$.
  By~\eqref{eq:cuspangle}, the relevant cusp angles are
  \begin{equation*}
    \measuredangle_{-\frac11}\biggl(0,-\frac12\biggr)
    =1\;,\qquad
    \measuredangle_{-\frac12}\biggl(-1,-\frac13\biggr)
    =2\;,\qquad\text{and}\qquad
    \measuredangle_{-\frac13}\biggl(-\frac12,\infty\biggr)
    =\frac23\;,
  \end{equation*}
  with sum~$\ell(P)=\frac{11}3$.
  Now, we can confirm the computation above because
  \begin{equation*}
    \bar\nu(M)
    =\frac{72}\pi\,A_{\mathrm{hyp}}(P)-24\,\ell(P)
    -\frac{72}\pi\,(\pi-2\thet)+3m_\rho(L;N_+,N_-)\\
    =-19\;.
  \end{equation*}
\end{ex}

\bibliographystyle{amsinitial}
\bibliography{g2geom}

\providecommand{\bysame}{\leavevmode\hbox to3em{\hrulefill}\thinspace}
\providecommand{\MR}{\relax\ifhmode\unskip\space\fi MR }
\providecommand{\MRhref}[2]{%
  \href{http://www.ams.org/mathscinet-getitem?mr=#1}{#2}
}
\providecommand{\href}[2]{#2}
\begin{thebibliography}{10}

\bibitem{atiyah75}
M.~F. Atiyah, V.~K. Patodi, and I.~M. Singer, \emph{Spectral asymmetry and
  {R}iemannian geometry, {I}}, Math. Proc. Camb. Phil. Soc. \textbf{77} (1975),
  97--118.

\bibitem{atiyah68}
M.~F. Atiyah and I.~M. Singer, \emph{The index of elliptic operators. {III}},
  Ann. of Math. (2) \textbf{87} (1968), 546--604.

\bibitem{beauville04}
A.~Beauville, \emph{Fano threefolds and {$K3$} surfaces}, The {F}ano
  {C}onference, Univ. Torino, Turin, 2004, pp.~175--184.

\bibitem{BCh4}
J.-M. Bismut and J.~Cheeger, \emph{Remarks on the index theorem for families of
  {D}irac operators on manifolds with boundary}, Differential geometry, Pitman
  Monogr. Surveys Pure Appl. Math., vol.~52, Longman Sci. Tech., Harlow, 1991,
  pp.~59--83.

\bibitem{BZtor}
J.-M. Bismut and W.~Zhang, \emph{An extension of a theorem by {C}heeger and
  {M}\"uller. {W}ith an appendix by {F}ran\c cois {L}audenbach}, Astérisque
  \textbf{205} (1992), 235 pp.

\bibitem{bunke95}
U.~Bunke, \emph{On the gluing problem for the $\eta$-invariant}, J. Diff. Geom.
  \textbf{41} (1995), 397--448.

\bibitem{BuMa}
U.~Bunke and X.~Ma, \emph{Index and secondary index theory for flat bundles
  with duality}, Aspects of boundary problems in analysis and geometry, Oper.
  Theory Adv. Appl., vol. 151, Birkhäuser, Basel, 2004, pp.~265--341.

\bibitem{cym}
A.~Corti, M.~Haskins, J.~Nordström, and T.~Pacini, \emph{Asymptotically
  cylindrical {C}alabi-{Y}au 3-folds from weak {F}ano 3-folds}, Geom. Topol.
  \textbf{17} (2013), 1955--2059.

\bibitem{g2m}
\bysame, \emph{$\textup{G}_2$-manifolds and associative submanifolds via
  semi-{F}ano 3-folds}, Duke Math. J. \textbf{164} (2015), 1971--2092.

\bibitem{eta}
D.~Crowley, S.~Goette, and J.~Nordström, \emph{An analytic invariant of
  ${G}_2$-manifolds}, arXiv:1505.02734v2, 2018.

\bibitem{nu}
D.~Crowley and J.~Nordström, \emph{New invariants of ${G}_2$-structures},
  Geom. Topol. \textbf{19} (2015), 2949--2992.

\bibitem{exotic}
\bysame, \emph{Exotic ${G}_2$-manifolds}, arXiv:1411.0656, 2018.

\bibitem{7class}
\bysame, \emph{The classification of 2-connected 7-manifolds}, Proc. Lond.
  Math. Soc., doi:10.1112/plms.12222, arXiv:1406.2226, 2019.

\bibitem{daifreed}
X.~Dai and D.~Freed, \emph{{APS} boundary conditions, eta invariants and
  adiabatic limits}, J. Math. Phys. \textbf{35} (2001), 5155--5194.

\bibitem{eells62}
J.~Eells, Jr. and N.~Kuiper, \emph{An invariant for certain smooth manifolds},
  Ann. Mat. Pura Appl. (4) \textbf{60} (1962), 93--110.

\bibitem{goette11}
S.~Goette, \emph{Adiabatic limits of {S}eifert fibrations, {D}edekind sums, and
  the diffeomorphism type of certain 7-manifolds}, J.~Eur.~Math.~Soc. (2014),
  2499--2555.

\bibitem{nuxx}
S.~Goette and J.~Nordström, \emph{$\nu$-invariants of extra twisted connected
  sums}, with an appendix by D. Zagier, in preparation, 2019.

\bibitem{gray70}
A.~Gray and P.~S. Green, \emph{Sphere transitive structures and the triality
  automorphism}, Pacific J. Math. \textbf{34} (1970), 83--96.

\bibitem{hhn}
M.~Haskins, H.-J. Hein, and J.~Nordström, \emph{Asymptotically cylindrical
  {C}alabi-{Y}au manifolds}, J.~Diff. Geom. \textbf{101} (2015), 213--265.

\bibitem{joyce00}
D.~Joyce, \emph{Compact manifolds with special holonomy}, OUP Mathematical
  Monographs series, Oxford University Press, 2000.

\bibitem{KiLe}
P.~Kirk and M.~Lesch, \emph{The $\eta$-invariant, {M}aslov index, and spectral
  flow for {D}irac-type operators on manifolds with boundary}, Forum Math.
  \textbf{16} (2004), 553--629.

\bibitem{kovalev03}
A.~Kovalev, \emph{Twisted connected sums and special {R}iemannian holonomy}, J.
  reine angew. Math. \textbf{565} (2003), 125--160.

\bibitem{mathai86}
V.~Mathai and D.~Quillen, \emph{Superconnections, {T}hom classes, and
  equivariant differential forms}, Topology \textbf{25} (1986), no.~1, 85--110.

\bibitem{milnor56}
J.~W. Milnor, \emph{On manifolds homeomorphic to the 7-sphere}, Ann. of Math.
  (2) \textbf{64} (1956), no.~2, 399--405.

\bibitem{milnor73}
J.~W. Milnor and D.~Husemöller, \emph{Symmetric bilinear forms},
  Springer-Verlag, New York, 1973, Ergebnisse der Mathematik und ihrer
  Grenzgebiete, Band 73.

\bibitem{milnor74}
J.~W. Milnor and J.~D. Stasheff, \emph{Characteristic classes}, Princeton
  University Press, Princeton, N. J., 1974, Annals of Mathematics Studies, No.
  76.

\bibitem{nikulin79}
V.~Nikulin, \emph{Integer symmetric bilinear forms and some of their
  applications}, Izv. Akad. Nauk SSSR Ser. Mat. \textbf{43} (1979), 111--177,
  238, English translation: {\it Math. USSR Izvestia} {\bf 14} (1980),
  103--167.

\bibitem{xtcs}
J.~Nordström, \emph{Extra-twisted connected sum ${G}_2$-manifolds},
  arXiv:1809.09083, 2018.

\bibitem{schelling}
A.~Schelling, \emph{{Die topologische $\eta$-Invariante und
  Mathai-Quillen-Str\"ome}}, Diploma thesis, Universit\"at Freiburg, {\tt
  http://www.freidok.uni-freiburg.de/volltexte/9530/}, 2014.

\bibitem{wallis18}
D.~Wallis, \emph{Disconnecting the moduli space of ${G}_2$-metrics via
  ${U}(4)$-coboundary defects}, arXiv:1808.09443, 2018.

\bibitem{wilkens71}
D.~L. Wilkens, \emph{Closed {$(s{-}1)$}--connected {$(2s{+}1)$}--manifolds},
  Ph.D. thesis, University of Liverpool, 1971.

\end{thebibliography}

\end{document}